\documentclass[11pt, a4paper]{article}

\usepackage{graphicx}
\usepackage{color}
\usepackage{amsmath}
\usepackage{amssymb}
\usepackage{amscd}
\usepackage{bbm}
\usepackage{tikz}
\usetikzlibrary{patterns}
\usepackage{hyperref}
\hypersetup{
    bookmarks=true,         
    colorlinks=true,       
    linkcolor=red,          
    citecolor=green,        
    filecolor=magenta,      
    urlcolor=cyan           
}

\newtheorem{Theorem}{Theorem}[section]
\newtheorem{Lemma}[Theorem]{Lemma}

\newtheorem{Definition}[Theorem]{Definition}

\newtheorem{Proposition}[Theorem]{Proposition}

\newtheorem{Remark}[Theorem]{Remark}

\newtheorem{Assumption}{Assumption}[section]

\numberwithin{equation}{section}

\newcommand{\inr}[1]{\bigl< #1 \bigr>}

\newcommand{\norm}[1]{\|#1\|}%

\newcommand{\beginproof}{{\bf Proof. {\hspace{0.2cm}}}}
\def \endproof
{{\mbox{}\nolinebreak\hfill\rule{2mm}{2mm}\par\medbreak}}

\DeclareMathOperator*{\argmin}{argmin}

\def\ds1{\textrm{1\kern-0.25emI}} 

\newcommand\E{{\mathbb E}}
\newcommand\R{{\mathbb R}}

\newcommand \cD{{\cal D}}

\newcommand \cF{{\cal F}}

\newcommand \cH{{\cal H}}

\newcommand \cL{{\cal L}}
\newcommand \cM{{\cal M}}
\newcommand \cN{{\cal N}}

\newcommand \cQ{{\cal Q}}
\newcommand \cR{{\cal R}}

\newcommand \cX{{\cal X}}
\newcommand \cY{{\cal Y}}

\newcommand \bP{{\mathbb P}}

\newcommand \bX{{\mathbb X}}

\usepackage{color}
\usepackage{ulem}
\newcommand{\im}{{\rm Im}}
\newcommand{\rk}{{\rm Rank}}

\begin{document}
\title{{Minimax regularization}}
\date{}
\author{Rapha{\"e}l Deswarte${}^{1,3}$ and Guillaume Lecu\'e${}^{2,4, 5}$}

 \footnotetext[1]{CMAP. Ecole Polytechnique, 91120 Palaiseau. Universit\'e Paris-Saclay. France}
 \footnotetext[2]{(corresponding author) CNRS, CREST, ENSAE, 5 avenue Henry Le Chatelier. 91120 Palaiseau. France.}
  \footnotetext[3] {Email: raphael.deswarte@polytechnique.edu}
 \footnotetext[4] {Email: guillaume.lecue@ensae.fr}
 \footnotetext[5]{Supported by the French National Research Agency (ANR) under the grant Labex Ecodec (ANR-11-LABEX-0047) and by the Maison des sciences de l'Homme under the grant 17MA01}

\maketitle

\begin{abstract}
Classical approach to regularization is to design norms enhancing smoothness or sparsity and then to use this norm or some power of this norm as a regularization function. The choice of the regularization function (for instance a power function) in terms of the norm is mostly dictated by computational purpose rather than theoretical considerations. 

In this work, we design regularization functions that are motivated by theoretical arguments. To that end we introduce a concept of optimal regularization called ``minimax regularization'' and, as a proof of concept,  we show how to construct such a regularization function for the $\ell_1^d$ norm for the random design setup. We develop a similar construction for the deterministic design setup.  It appears that the resulting regularized procedures are different from the one used in the LASSO in both setups.
\end{abstract}

\textbf{Keywords:} Regularization, Gaussian mean width, empirical processes, regularization.


\section{Introduction} 
\label{sec:introduction}


Let $(\cX, \mu)$ be a probability space and $(X,Y)$ be a couple of random variables, in which $X$ is distributed according to $\mu$. One is given a sample of $N$ independent couples $(X_i,Y_i)_{i=1..N}$ distributed according to the joint law of $(X,Y)$. On the basis of this sample, one tries to link $X$ and $Y$ by a random mapping $\hat{f}$ with $\hat{f}(X)$ close (in $L_2$) to $Y$. This is the classical problem, in learning theory, of the prediction of an output $Y$ from an input $X$ given a i.i.d. copies of the couple $(X,Y)$. 

To that end, one is given a class $F$ of functions from $\cX$ to $\R$ and the aim in learning theory is to mimic the \textit{best} element in $F$ for the prediction of $Y$ by a function of $X$ in $F$. We assume that $F$ is  closed and convex in $L_2(\mu)$ so that it exists a function $f^*$ that minimizes the square loss in $F$:
\begin{equation}
\label{eq:oracle}
f^*\in\argmin_{f\in F}\E \left(Y-f(X)\right)^2.
\end{equation}This function is usually called the oracle (cf. \cite{MR1775638}); it is the closest function in $F$ to $Y$ in $L_2$. Now, the goal is to construct an estimator $\hat{f}$ whose $L_2(\mu)$ distance to $f^*$ is as small as possible  using the dataset $\{(X_i, Y_i):i=1, \ldots, N\}$. In the framework considered in this paper, the excess risk of $\hat f$, which is the difference $\E(Y-\hat f(X))^2-\E (Y-f^*(X))^2$, is actually equal to $\norm{\hat f - f^*}_{L_2(\mu)}^2$ and so estimating $f^*$ is equivalent to predicting $Y$; thus we fall back on the original prediction problem by estimating $f^*$ in $L_2(\mu)$.

One may therefore try to bound the quadratic error $\norm{\hat f - f^*}_{L_2(\mu)}$ either \textit{in expectation} or \textit{in deviation} with respect to the sample. In this work, we obtain upper bounds on the quadratic error that are valid in deviation, showing that the results are true for ``most" samples rather than in average.

Given that we want to be close to a function $f^*$ minimizing $f\to \E(Y-f(X))^2$ over $F$, a natural candidate for this problem is the Empirical Risk Minimizer (ERM) also known as  the ``least squares estimator": 
\begin{equation}
\label{eq:erm}
\hat{f}_{ERM} \in\argmin_{f\in F}\frac{1}{N}\sum_{i=1}^N(Y_i-f(X_i))^2.
\end{equation}
 Many works have been carried out for general classes (see, \cite{MR2829871,MR2319879,MR1739079,MR1385671}) or on the vectorial case(see \cite{stein1956inadmissibility} 
 for the famous Stein paradox, and \cite{chatterjee2014new} for elements about the  admissibility of the ERM).

 It appears that when F is too large (for instance the whole $L_2(\mu)$ space), the ERM tends to ``overfit". The understanding of this phenomenon has led to the introduction of ``regularization methods" which were originally used to smooth estimators in order to overcome the ``overfitting phenomena". Those procedures are nowadays used beyond their smoothing effect and, in particular,  they are now extensively used in Statistics and learning theory for their ``low-dimensional / sparsity inducing properties''. At a high level description, those methods make a trade-of between an ``adequation to the data term'' and a ``regularization term'' and their general form (for the quadratic loss) is 
\begin{equation}\label{eq:RERM}
\hat f \in\argmin_{f\in F}\Big(\frac{1}{N}\sum_{i=1}^N(Y_i-f(X_i))^2 + \Psi(f)\Big)
\end{equation}where $\Psi$ is a function usually called the \textit{regularization function}.

The ``adequation to the data term'' can be constructed from any loss function; for the case of the quadratic loss, this term reads like  $N^{-1}\sum_{i=1}^N(Y_i-f(X_i))^2$.

As for the regularization term $\Psi$, several choices are possible, enabling to smooth the estimator, or to force a low-dimensional structure. It depends thus on the a priori knowledge one has on the data (and in particular on $f^*$), and on computational issues.

A first option is the Tikhonov / ridge regularization: 
\begin{equation}\label{eq:tikhonov}
 \hat f_\lambda \in \argmin_{f\in\cH}\Big(\frac{1}{N}\sum_{i=1}^N (Y_i-f(X_i))^2 + \lambda \norm{f}_\cH^2\Big)
 \end{equation} where $\norm{\cdot}_{\cH}$ is a Hilbert norm. One expects in this case that $\norm{f}_\cH$ reflects the smoothness of $f$ (for instance $\norm{f}_\cH = f(0) + \big(\int_{\R}|f^\prime(t)|^2dt\big)^{1/2}$) and that the oracle function $f^*$ 
has a small $\norm{f^*}_\cH$ norm.

In the finite but high dimensional vectorial case, one often wishes the estimator to be sparse, i.e., to have few non-zero components (for some well designed basis, cf. \cite{MR2479996}). This may come from the fact that the vector to be estimated is known in advance to be sparse; or that, in high-dimensional problems,  it is computationally important not to have to manage a huge amount of non-zero coefficients. 

Then, a natural way to address this question is to use a sparsity-inducing penalization: like the number of non-zero components of the vector, sometimes called the ``$\ell_0$- norm''. Even though it is theoretically appealing (cf. \cite{giraud2014introduction}), it proves to be computationally intractable (actually NP-Hard, in general, cf. \cite{natarajan1995sparse}). But for geometric reasons, another regularization is efficient to induce sparsity: the $\ell_1$ norm (which can be see as the convex relaxation of the $\ell_0$ norm on the unit $\ell_\infty$-ball). The associated estimator is called the LASSO (``Least Absolute Shrinkage and Selection Operator", \cite{MR1379242}):
 \begin{equation}\label{eq:lasso}
 \hat t_\lambda \in \argmin_{t\in\R^d}\Big(\frac{1}{N}\sum_{i=1}^N (Y_i-\inr{X_i,t})^2 + \lambda \norm{t}_1\Big).
 \end{equation}
In \cite{candes2009near}, it is emphasized that Lasso leads generally to sparsity, though giving some counter-examples in which it doesn't work well --in particular Lasso struggles when dealing with a data matrix with high correlations among the columns. To tackle this kind of issues, it is possible to ``mix" regularizations: it is the principle of the ``Elastic net method" (\cite{MR2137327}), which  penalty is a combination of the $\ell_1$ and $\ell_2$ norms.

All these methods rely on the choice of one (even two for the Elastic net) regularization parameter $\lambda$, fixed by the  statistician on the basis of empirical methods such as cross-validation. It has to be chosen wisely in order to make the right trade-off between the two terms and thus to minimize the rate of convergence of the regularization procedure towards the oracle.
\subsection{Regularization and Model Selection}
At this point, it should be clear to the reader that choosing (or even designing) a specific regularization norm like $\norm{\cdot}_\cH$ or $\norm{\cdot}_1$ depends only on an a priori knowledge we have. But once this choice has been made, why would someone use the square of this norm in one case (as for the Tikhonov / ridge regularization), or the norm itself (as for the LASSO) or some other power of this norm (cf. \cite{MR2816342} for some examples) in other situations? In many cases, this choice is only made following some computational considerations. 

In this work, we want to support the choice of regularization functions (given some norm) on theoretical arguments. To that end we will rely on the model selection theory and, in particular, on a key principle in model selection which is to design penalty functions that capture the ``complexity'' of a model in the most accurate way.  

The right calibration of penalty functions has opened an important stream of researches since the work \cite{MR1679028}. It has led many researchers to (re)think about the notion of complexity in statistics. In a nutshell, there are mainly three types of quantities that have been introduced to measure the statistical complexity of a statistical model: combinatorial (like the VC dimension, \cite{MR1641250}), metric (like the entropy) and random (like Gaussian mean width and Rademacher complexities, \cite{MR2240689,MR2329442}). Penalty calibration has culminated with the notion of ``minimal penalty'' that are sharp penalty functions with exact constants (cf. for instance \cite{birge2001gaussian}) thanks to second order term analysis of the notion of complexity of a model.

In the present work, we want to put forward the idea that the ``right'' choice (from a theoretical point of view) of a regularization function may also follow from a careful study of the complexity of a specific family of models. To that end, the first argument is to look at regularization as a model selection problem for which one has to design a sharp penalty. This has been done for instance in Chapter~3.7 in \cite{HDR}.
 In the particular case where one is given a norm $\norm{\cdot}$  for regularization, then the associated regularized ERM is a penalized estimator associated to the sequence of embedded models $(m_r)_{r\geq0}$ where for all $r\geq0$, $m_r=\{f\in F:\norm{f}\leq r\}$ and the right way to regularize is given by a function ${\rm reg}: f\in F\to {\rm pen}(m_{\norm{f}})$ where $m_{\norm{f}}$ is the smallest model in $(m_r)_{r\geq0}$ containing $f$, see \cite{birge2001gaussian}. This idea is a baseline of this work. 

Before diving into further details about the way we suggest to construct a  regularization function, let us precise what we expect from a good procedure, in particular how we evaluate that a regularization function is the ``right" one, at least from a theoretical point of view.   We therefore need to introduce a concept of optimality for regularized estimators. Once again we rely on the basics of model selection theory.

Model selection procedures have been used originally to construct adaptive estimators. For the model selection problem we want to solve, this adaptivity problem reads like selecting the smallest model in the family $\left(\{f\in F:\norm{f}\leq r\}\right)_{r\geq0}$ containing $f^*$ which is obviously $\{f\in F:\norm{f}\leq \norm{f^*}\}$. Therefore, the adaptation problem we want to solve here is to construct a procedure which performance is as good as if we had been given the value $\norm{f^*}$ in advance. In particular, an estimator achieving the minimax rate of convergence over the model $\{f\in F:\norm{f}\leq \norm{f^*}\}$ would solve this adaptation problem. In what follows, we design regularization functions in order to meet this requirement but before that we clarify the notion of minimax rate over a model for the type of deviation results we prove below. 

To simplify the exposition, we will focus on a specific, though very classical and widely-used, framework:  the vectorial case, i.e. when $F=\left\{\inr{\cdot, t}: t\in T\right\}$ is a class of linear functionals from $\R^d$ to $\R$ indexed by some subset $T\subset\R^d$, with Gaussian design, and Gaussian noise (with known variance $\sigma^2$).

\begin{Definition}\label{def:minimax}
Let $T\subset\R^d$,  $X$ denote a standard Gaussian vector in $\R^d$ and $\xi$ be a centered real-valued Gaussian random variable  with variance $\sigma^2$, independent of $X$. For all $t^*\in T$, define the random variable $Y^{t^*} = \inr{X, t^*} + \xi$ and denote by $\cY^T:=\{Y^{t^*}:t^*\in T\}$ the set of all such random variables. 

Let $\hat t_N$ be a statistics from $(\R^d\times \R)^N$ to $\R^d$. Let $0<\delta_N<1$  and $\zeta_N>0$. We say that \textbf{$\hat t_N$ performs with accuracy $\zeta_N$ and confidence $1-\delta_N$ relative}\textbf{ to the set of targets $\cY^T$}, if for all $Y\in\cY^T$, with probability, w.r.t. to a sample $\cD:=\{(X_i, Y_i): i=1,\cdots, N\}$ of i.i.d. copies of $(X, Y)$, at least $1-\delta_N$,  $\norm{\hat t_N - t^*}_2^2\leq  \zeta_N$.

We say that \textbf{$\cR_N$ is a minimax rate of convergence over $T$ for the confidence $1-\delta_N$} if the two following statements hold:
\begin{enumerate}
	\item there exists a statistics $\hat t_N$ which performs with accuracy $\cR_N$ and confidence $1-\delta_N$ relative to the set of targets $\cY^T$
	\item there exists an absolute constant $g_0>0$ such that if $\tilde t_N$ is a statistics which attains an accuracy $\zeta_N$ with confidence $1-\delta_N$ relative to the set of targets $\cY^T$ then necessarily $\zeta_N \geq g_0 \cR_N$.
\end{enumerate}
\end{Definition}


In the following (cf. Theorem~\ref{thm:minimax_B1} below), we recall a result from \cite{LM13} on the minimax rate of convergence over $T= \rho B_1^d$, the unit $\ell_1^d$-ball of radius $\rho\geq0$, for a constant confidence (i.e., for instance,  $\delta_N=1/4$).  Note that classical minimax rates of convergence are usually given in expectation (cf. for instance \cite{MR2724359}). The main difference here with Definition~\ref{def:minimax} is that it is given for deviation results: the minimax rate $\cR_N$ may depend on the confidence parameter $\delta_N$ (cf. \cite{LM13}).     

In the present work, we are interested in procedures achieving the minimax rate of convergence over the model $\{t\in \R^d: \norm{t}\leq \norm{t^*}\}$. This provides a natural way to introduce a notion of optimality for the problem of designing regularization functions.


\begin{Definition}\label{def:optimality}
Let $\norm{\cdot}$ be a norm on $\R^d$, $0<\delta_N<1$ and $T\subset\R^d$.   Let us consider the following RERM for some function $\Psi:\R_+\to\R$: 
\begin{equation*}
 \hat t \in\argmin_{t\in \R^d}\left(\frac{1}{N}\sum_{i=1}^N (Y_i-\inr{X_i, t})^2+\Psi(\norm{t})\right)
\end{equation*}constructed from a sample $\cD:=\{(X_i, Y_i): i=1,\cdots, N\}$ of i.i.d. copies of $(X, Y^{t^*})$ where $Y^{t^*}=\inr{X,t^*} + \xi$ with $X\sim\cN(0,I_{d\times d})$, $\xi\sim\cN(0, \sigma^2)$ is independent of $X$ and $t^*\in \R^d$. We say that $\Psi$ is a \textbf{minimax regularization function for the norm $\norm{\cdot}$ and the confidence $1-\delta_N$ over $T$}, if there exists an absolute constant $g_1>0$ such that for all $t^*\in T$, the RERM $\hat t$ is such that with $\bP_{t^*}$-probability at least $1-\delta_N$, $\norm{\hat t - t^*}_2^2 \leq g_1\cR_N^{\norm{t^*}_1}$, where $\cR_N^{\norm{t^*}_1}$ is the minimax rate of convergence over $\{t\in \R^d: \norm{t}\leq \norm{t^*}\}$ and $\bP_{t^*}$ denotes the probability distribution of a $N$ sample of i.i.d. copies of $(X, Y^{t^*})$.
\end{Definition}

The aim of this work is to show that one can design minimax regularization functions by finding the right notion of complexity of the sequence of embedded models $\left(\{t\in\R^d:\norm{t}\leq r\}\right)_{r\geq0}$. Note however that there should be some situations where designing such an optimal regularization function would be impossible at some given confidence parameter $\delta_N$. In particular, such a situation should happen when the  Empirical risk minimization (ERM) procedure over the ``true model'' $\{t\in \R^d:\norm{t}\leq \norm{t^*}\}$ is not itself a minimax procedure over the model $\{t\in \R^d:\norm{t}\leq \norm{t^*}\}$. This happens for constant confidence bound (for instance, when $\delta_N=1/4$) when there is a gap in Sudakov inequality (cf. \cite{LM13} for more details). Nevertheless, in \cite{LM13}, it is proved that for high confidence bounds (that is when $\delta_N$ decays exponentially fast with the complexity of the model) ERM is always minimax over convex classes. 
It appears that for the case of $\ell_1^d$-balls ERM is minimax for all confidence regime therefore this subtlety will not show up in this special case.

\subsection{General approach provided in this paper}

Let us now present our approach.  As we mentioned before, we want to construct a regularization function depending on the complexity of the models $\{t\in\R^d: \norm{t}\leq r\}$ for all $r\geq0$. This leads to choose  a RERM (in the vectorial case) having the following form:

\begin{equation}\label{eq:model_selection}
\hat t \in\argmin_{t\in\R^d}\Big(\frac{1}{N}\sum_{i=1}^N(Y_i-\inr{X_i, t})^2 + {\rm comp} \left(\norm{t}B_{\norm{\cdot}}\right)\Big)
\end{equation}where $B_{\norm{\cdot}}$ is the unit ball associated with the given regularization norm $\norm{\cdot}$ and for all $t\in\R^d$, $\norm{t}B_{\norm{\cdot}}=\{u\in\R^d: \norm{u}\leq \norm{t}\}$. The key feature in \eqref{eq:model_selection} is the ``complexity function'' $r\geq 0 \to {\rm comp} \left(rB_{\norm{\cdot}}\right)$ which aims at measuring with the best possible accuracy the complexity of the models $rB_{\norm{\cdot}}$  for all $r\geq0$ from a statistical point if view. Or course, finding the right notion of complexity is paramount in this approach. 

To aim at optimality, in this paper we advocate complexities that are tailored for the specific statistical problem at stake. It turns out that the ``right" choice of complexity, and hence, of regularization, is linked to the behavior of two empirical processes. Those two empirical processes are ultimately connected to the two sources of statistical complexities in the considered problem. When estimating $t^*$ from the data $(X_i, Y_i)_{i=1}^N$ there are two statistical issues : 1)(an inverse problem) $t^*$ is observed only through $\bX$, where $\bX\in\R^{N\times d}$ is the operator whose rows vectors are given by the $X_i$'s; 2)(noisy data) the observations have been corrupted by some noise $\xi$. The action of the operator $\bX$ on the models $r B_{\norm{\cdot}}$ for all $r\geq0$ plays a prominent role in our analysis. In particular, the size of the intersection of its kernel with the model is a natural minmax lower bound for any estimator since any two vectors in the kernel of $\bX$ and the model are indistinguishable. The effect of the ``distortion'' of the operator $\bX$ does not show up for small models (i.e. small values of $r$) because of the presence of the noise which blurs everything at small scales. But passing beyond some threshold for the signal-to-noise ratio $r/\sigma$, only the distortion of $\bX$ matters from a statistical point of view. This phenomenon occurs only when $N\lesssim d $, because that is the regime where $\bX$ has a non trivial kernel. On the contrary, in the low dimensional setup $d\lesssim N$, $\bX$ is well conditioned with high probability and therefore, there is no statistical complexity coming from the distortion of $\bX$ since in that regime there is no such distortion. Controlling the distortion of $\bX$ is a key issue in high-dimensional statistics. It is behind all classical properties like RIP (cf. \cite{MR2723472}) or REC (cf. \cite{MR2533469}) and it will play equally a key role in our analysis. In particular, the $\ell^d_2$ diameter of the intersection of $\bX$ with the model $r B_{\norm{\cdot}}$ will appear explicitly in the optimal regularization. Given that $\bX$ is a standard Gaussian random matrix, this diameter will be the Gelfand width of $r B_{\norm{\cdot}}$ in our example (cf. \cite{MR774404}, Chapter~2 in \cite{MR3113826} or \cite{LM13} for more details on Gelfand widths and their role in signal processing and learning theory).  


\subsection{Overview of the paper; main results}

As a proof of concept we present an example of the construction of a minimax regularization function in the popular set-up of regularization by the $\ell_1^d$ norm. Let us recall the  statistical model we used in both Definition~\ref{def:optimality} and Definition~\ref{def:minimax}: the Gaussian linear regression model with a  Gaussian design
\begin{equation}\label{eq:model}
Y=\inr{X,t^*} + \xi
\end{equation}where $X\sim \cN(0, I_{d\times d})$ and $\xi\sim\cN(0,\sigma^2)$ are independent, centered Gaussian variables in $\R^d$ and $\R$ respectively. As written before, a dataset $(X_i,Y_i)_{i=1}^N$ of $N$ i.i.d. copies of the couple $(X,Y)$ is provided and one wants to use it to estimate $t^*$.


 Note that we choose a Gaussian random design to make the exposition as simple as possible. The results can be extended to more general sub-Gaussian designs. Nonetheless, our goal is not to provide general results but to show that the approach we present allows to achieve minimax regularization in some classical set-up. Moreover, we want to see the effect of the random design on the construction of a minimax regularization.  

In the supplementary material~\ref{fixed_design}, we consider a fixed design setup, in which one still has for all $i$: $Y_i=\inr{X_i,t^*} + \xi_i$ with $\xi_i\sim\cN(0,\sigma^2)$ i.i.d., but with $X_i$ deterministic, satisfying 
an ``isomorphic property'' on ``compressible vectors'', equivalent to the RIP from \cite{MR2723472}. We will see that under this property, the arguments and the results will be quite similar to the random design case.

We will not be interested in getting optimal or sharp numerical constants, and some of the inequalities and coefficients in the arguments will be rather loose from this point of view, but what actually matters is that the quantity will have the right order of magnitude w.r.t. $N, d, \sigma$ and $\norm{t^*}_1$. 

As for our choice to consider regularizations that are functions of the $\ell_1^d$ norm, its motivation is that the $\ell_1^d$-norm has been one of the most studied regularization norm since the beginning of high-dimensional statistics, in particular for the reasons presented in the Introduction. Moreover, as mentioned previously, the ERM over $\ell_1^d$-balls is minimax for every confidence $1-\delta_N$ (cf. \cite{LM13}), this makes the construction of minimax regularization possible for different deviation parameters and this makes the exposition also simpler.

 Let the choice of the $\ell_1^d$-norm as a regularization norm be made once and for all. Now, the problem we want to solve is to construct a regularization function $\Psi:\R\to\R$ such that the regularized procedure
\begin{equation}\label{eq:minimax_lasso}
 \hat t_\Psi \in \argmin_{t\in \R^d}\Big(\frac{1}{N}\sum_{i=1}^N (Y_i-\inr{X_i,t})^2 + \Psi(\norm{t}_1)\Big)
 \end{equation}achieves the minimax rate of convergence over $\norm{t^*}_1 B_1^d$ given $N$ i.i.d. data $(X_i, Y_i), i=1, \ldots, N$ distributed according to \eqref{eq:model}. And,  we want $\hat t_\Psi$ to satisfy that property whatever $t^*\in\R^d$ is.
 
Now, let us explain the strategy we use to design a minimax regularization function. We denote for all $\rho\geq0$ and $r\geq0$,  $\rho B_1^d = \{t\in\R^d: \norm{t}_1\leq \rho\}$, and  $rB_2^{d} = \{t\in\R^d: \norm{t}_2\leq r\}$. The starting point to our approach is that $\hat t_\Psi$ minimizes $t\mapsto P_N\cL_t^\Psi$ over $\R^d$, where, for every $t\in \R^d$, 
$$P_N\cL_t^\Psi:=\left(\frac{1}{N}\sum_{i=1}^N (Y_i-\inr{X_i,t})^2 + \Psi(\norm{t}_1)\right)- \left(\frac{1}{N}\sum_{i=1}^N (Y_i-\inr{X_i,t^*})^2 + \Psi(\norm{t^*}_1)\right),$$  in particular, $P_N\cL_{\hat t}^\Psi \leq P_N\cL_{t^*}^\Psi=0$. So if one shows that $P_N \cL_t^\Psi>0$ for all $\norm{t}_1\gtrsim \norm{t^*}_1$ then this will prove that $\norm{\hat t_\Psi}_1\lesssim \norm{t^*}_1$, proving that $\hat t_\Psi$ belongs to the right model. This will be essentially the main step since for the correct choice of $\Psi$, we will show that the regularization has no effect within the right model and that the RERM $\hat t^\Psi$ has essentially the same statistical behavior as the ERM in $\norm{t^*}_1 B_1^d$ which is known to be minimax. We will therefore conclude that $\hat t^\Psi$ can learn $t^*$ at the minimax rate of convergence within the model $\norm{t^*}_1 B_1^d$ without knowing in advance the radius $\norm{t^*}_1$.

 Using the quadratic / multiplier decomposition as in \cite{LM13,saum:12,saumard2017new}, one can write $P_N\cL_t^\Psi$ as the sum of three terms: $P_N\cL_t^\Psi=P_N {\cal Q}_{t-t^*}+P_N {\cal M}_{t-t^*}+{\cal R}_{t,t^*}$ where
\begin{itemize}
\item $P_N {\cal Q}_{t-t^*}:=\sum_{i=1}^N \left(\inr{ X_i,t^*}-\inr{X_i,t}\right)^2/N$ is the ``quadratic process"
\item $P_N {\cal M}_{t-t^*}:=2\sum_{i=1}^N \left(Y_i-\inr{X_i,t^*}\right)\left(\inr{X_i,t^*}-\inr{X_i,t}\right)/N$ is  the ``multiplier process" 
\item ${\cal R}_{t,t^*}:=\Psi(\norm{t}_1)-\Psi(\norm{t^*}_1)$ is the regularization part.
\end{itemize}

The definition of our complexity will thus be a consequence of the study of the behavior of the quadratic and multiplier empirical processes indexed by $t\in \rho B_1^d$ for all $\rho\geq0$. The two processes are associated to the two statistical complexities previously discussed: 1) the quadratic process can be written as $P_N {\cal Q}_{t-t^*} = \norm{\bX(t-t^*)}_2^2$ and is well behaved (i.e. of the order of $\norm{t-t^*}_2^2$) when $\bX$ is well conditioned; 2) the multiplier process is measuring the statistical complexity coming from the noise $\xi=Y-\inr{X, t^*}$, $P_N \cM_{t-t^*}$ is the empirical correlation between the noise and the model shifted by $\inr{\cdot, t^*}$. All the game is now to identify regions of the space $\R^d$ where the statistical complexity come from the distortion of $\bX$ or from the noise. This drives the construction of the optimal regularization function $\Psi$.

 In order to identify those regions, note that for every fixed $t\in\R^d$, the distribution of these two processes depend on $t-t^*$ only by its $\ell_2^d$-norm $\norm{t-t^*}_2$,  in two different ways: $P_N {\cal Q}_{t-t^*}$ in a quadratic way, $P_N {\cal M}_{t-t^*}$ in a linear way. So it is natural to partition the model $\rho B_1^d$ into vectors with ``small" $\ell_2$ norm -- i.e. the intersection of $\rho B_1^d\cap r B_2^d$ for an adequate radius $r$ -- and vectors of $\rho B_1^d$ with $\ell_2^d$-norm larger than $r$. We will see that outside $rB_2^d$, with high probability the two processes are ``well-behaved" and regularization is unnecessary; but inside $rB_2^d$ it is not the case, the operator $\bX$ may have a kernel and the noise is making the estimation hard: hence, this is where the regularization will be needed to keep control of the situation and this is precisely where the regularization function is designed. In that case, either the statistical complexity comes from the size of the intersection of the kernel of $\bX$ with $\rho B_1^d$ and therefore one needs to take $\Psi(\rho)$ of the order of this diameter (which appears to be equal to the Gelfand width of $\rho B_1^d$ to the square) or the statistical complexity comes from the noise and then $\Psi(\rho)$ is of the order of the oscillations of the multiplier process inside $\rho B_1^d\cap r B_2^d$. 

The choice of the ``adequate radius" $r$ is of course paramount in our approach. It results from the right understanding of the two previously discussed sources of statistical complexities:  the bigger these complexities, the bigger this radius (since, as we mentioned, outside $rB_2^d$ the processes are well-behaved). First, we want to identify  the smallest $\ell_2^d$ radius $r_Q(\rho)$ above which $\bX$ is well-behaved in $\rho B_1^d$, i.e. such that for every $t\in\rho B_1^d$, if $\norm{t-t^*}_2\geq r_Q(\rho)$, then  $P_N \cQ_{t-t^*} = \norm{\bX(t-t^*)}_2^2\sim\norm{t-t^*}_2^2$. Then, we need to identify the smallest $\ell_2^d$-radius $r_M(\rho)$ above which the effect of the noise is below the signal intensity that is above which one can clearly identify if $t\neq t^*$ when $\norm{t-t^*}_2\geq r_M(\rho)$. To that end we want to make the oscillations of the multiplier process smaller than the one of the quadratic process, which is of the order of $\norm{t-t^*}_2^2$ when $\norm{t-t^*}_2\geq r_Q(\rho)$. It will appear that, in our framework, the two radii obtained from the above trade-offs are solution of fixed point equations for all $\rho\geq0$: for some absolute constants $Q$ and $\eta$ (to be chosen later): 
\begin{itemize}
\item the ``quadratic fixed point" is $r_Q(\rho):= \inf\Big(r>0:\ell^*(\rho B_1^d\cap r B_2^{d})= Q r \sqrt{N}\Big)$
\item the ``multiplier fixed point" is $r_M(\rho):= \inf\Big(r>0:\sigma \ell^*( \rho B_1^d\cap r B_2^{d})= \eta r^2 \sqrt{N}\Big)$
\end{itemize}
where $\ell^*( \rho B_1^d\cap rB_2^{d})$ is the Gaussian mean width of the localized set  $\rho B_1^d\cap r B_2^{d}$ defined as 
\begin{equation*}
\ell^*( \rho B_1^d\cap r B_2^{d}) = \E\sup_{t\in \rho B_1^d\cap r B_2^{d}}\inr{G, t}
\end{equation*} where $G$ is a standard Gaussian vector in $\R^d$. As our framework involves ``Gaussian randomness" in both the design and the noise, it is not surprising that the Gaussian mean width arise when dealing with the control of the two processes. However, Gaussian mean widths appear in learning theory, statistics and signal processing way beyond the ``full Gaussian framework'' as considered here (see, for instance, \cite{MR3612870}).

These two fixed points have been introduced in \cite{LM13} and used later in \cite{lecue2016regularization} for the study or ERM and RERM. As their names suggest, the quadratic fixed point will be used to control the quadratic process, and the multiplier fixed point to control the multiplier process. Their general definitions use an inequality rather than an equality inside the infimum: $r_Q(\rho)$ is defined as $\inf\Big(r>0:\ell^*(\rho B_1^d\cap r B_2^{d})\leq Q r \sqrt{N}\Big)$ and $r_M(\rho)$ as $\inf\Big(r>0:\sigma \ell^*( \rho B_1^d\cap r B_2^{d})\leq \eta r^2 \sqrt{N}\Big)$. This allows to deal with infinite-dimensional set-ups in which the mapping $r \mapsto \ell^*(\cF \cap r B_2) $ is not necessarily continuous. But in our case, this mapping is continuous and the infimum is attained in a point for which there is exact equality.

 It appears that one can provide an explicit formulation for the two fixed point $r_Q(\rho)$ and $r_M(\rho)$ in many situations and, in particular, in the case of the $\ell_1^d$-norm (cf. \cite{LM13}): for some absolute constants $C_M^{(1)}, C_M^{(2)}, C_Q^{(1)}, C_Q^{(2)}$ and $\zeta < 1 < \zeta^\prime$, for all $\rho$, there exists $C_M\in [C_M^{(1)}, C_M^{(2)}]$ such that:
\begin{equation}\label{eq:rM}
r_M^2(\rho) = C_M  \left\{
\begin{array}{cl}
\frac{\sigma^2 d}{N} & \mbox{ if } \rho^2 N \geq \sigma^2 d^2\\
\\
\rho \sigma\sqrt{\frac{1}{N}\log\Big(\frac{e\sigma d}{\rho\sqrt{N}}\Big)} & \mbox{ if } \sigma^2 \log d \leq \rho^2N\leq  \sigma^2 d^2\\
\\
\rho \sigma \sqrt{\frac{\log(ed)}{N}} & \mbox{ if } \rho^2N\leq \sigma^2 \log d.
\end{array}
\right.
\end{equation}
and there exists $C_Q\in[C_Q^{(1)},C_Q^{(2)}]$ such that
\begin{equation}\label{eq:rQ}
r_Q^2(\rho) = C_Q \left\{
\begin{array}{cl}
 0 & \mbox{ if }  N \geq \zeta^\prime  d \\
\frac{\rho^2}{N}\log\Big(\frac{ed}{N}\Big) & \mbox{ if } N\leq \zeta d.
\end{array}
\right.
\end{equation}Note that when $\zeta d \leq N\leq \zeta^\prime d$, $r_Q(\rho)$ decays from $(\rho^2/N)\log(ed/N)$ to $0$ and one only has an upper estimate on $r_Q(\rho)$ given by $C_Q\rho^2/N$.  We will therefore not consider this case in the following since it involves to deal with sharp estimates on the spectra of squared or approximatively squared Gaussian random matrices. Note also that $C_M$ and $C_Q$ may depend on $\rho, N, d$, $\sigma$  but they are both controlled from above and below by absolute constants (independent of $\rho, N, d$ and $\sigma$).

Now that we have a way to measure the statistical complexity of a model we need one more thing before turning to the effective construction of a minimax regularization for the $\ell_1^d$-norm: we need to know the minimax rate of convergence over $\ell_1^d$-ball $\rho B_1^d$ for all $\rho\geq0$. We will see below that one way to measure the statistical complexity of a model is closely related to its minimax rate. To that end, we summarize the main results in the constant deviation case $\delta_N=1/4$ from section~4.1 in ~\cite{LM13} in the following theorem. 

\begin{Theorem}\label{thm:minimax_B1}Consider the Gaussian linear model with Gaussian design introduced in \eqref{eq:model}.  Let $\rho>0$. The minimax rate of convergence for constant confidence parameter $\delta_N=1/4$ over $\rho B_1^d$ is achieved by the ERM and is given (up to absolute constants) by
\begin{equation}\label{eq:minimax_B1}
\min\left(r^2(\rho), \rho^2\right) \mbox{ where } r(\rho) = \max\big(r_Q(\rho), r_M(\rho)\big).
\end{equation}Up to multiplicative absolute constants, this rate is given  for some $\zeta <1<\zeta^\prime$,
\begin{enumerate}
  \item when $N\leq \log d$, by $\rho^2$,
  \item when $\log d\leq N\leq \zeta d$, by  
  \begin{equation*}
\begin{cases}
\rho^2 & \mbox{ if } \rho^2 N \leq \sigma^2 \log d,\\
\rho \sigma \sqrt{\frac{1}{N}\log\Big(\frac{ed^2\sigma^2}{\rho^2N}\Big)} & \mbox{ if }  \sigma^2 \log d\leq \rho^2N \leq \frac{\sigma^2 N^2}{\log(ed/N)} ,
\\
\frac{\rho^2}{N}\log\Big(\frac{ed}{N}\Big) & \mbox{ if } \rho^2N \geq \frac{\sigma^2 N^2}{\log(ed/N)}
\end{cases}
\end{equation*}
\item  when $N\geq \zeta^\prime d$, by
\begin{equation*}
\begin{cases}
\rho^2 & \mbox{ if } \rho^2 N \leq \sigma^2 \log d,\\
\rho \sigma \sqrt{\frac{1}{N}\log\Big(\frac{ed^2\sigma^2}{\rho^2N}\Big)} & \mbox{ if }  \sigma^2 \log d\leq \rho^2N \leq \sigma^2 d^2,
\\
\frac{\sigma^2 d}{N} & \mbox{ if } \rho^2 N \geq \sigma^2 d^2.
\end{cases}
\end{equation*}
\end{enumerate}In other words, for all $\rho\geq0$ and $t^*\in\rho B_1^d$,  the ERM $\hat t^{ERM}_\rho\in\argmin_{t\in\rho B_1^d} \sum_{i=1}^N (Y_i-\inr{X_i, t})^2$, is such that, with probability at least $3/4$, $\norm{\hat t^{ERM}_\rho - t^*}_2^2\leq \min\left(r^2(\rho), \rho^2\right)$. Moreover, there are no estimator that can do uniformly better than the ERM $\hat t^{ERM}_\rho$ over $\rho B_1^d$ when $N\notin(\zeta d, \zeta^\prime d)$.
\end{Theorem}

Note that we have decided to present the result in the constant deviation result (that is for $\delta_N=1/4$) whereas it is actually true with a much better probability estimate in section~4.1 in ~\cite{LM13}. We will also obtain our main results with an exponentially large deviation below.

As mentioned previously, when $N\in[\zeta d, \zeta^\prime d]$, we only have an upper bound on $(r_Q(\rho))^2$ that does not match the minimax lower bound. As a consequence, the $N\sim d$ regime is not considered in Theorem~\ref{thm:minimax_B1}. Notable is that the rate $\rho^2$ is the trivial rate obtained by taking the $\ell_2^d$ diameter of the model $\rho B_1^d$ which is simply $2\rho$. Therefore, any statistics $\tilde t_N$ (like the ERM $\hat t^{ERM}_\rho$) taking its values in $\rho B_1^d$ satisfies with probability $1$, $\norm{\tilde t_N - t^*}_2^2\leq 4 \rho^2$ for all $t^*\in \rho B_1^d$. This is a  trivial bound that one can get for free as long as the radius $\rho$ is known. However, for the construction of an optimal regularization function which can be seen as an adaptation to the radius $\norm{t^*}_1$, which is therefore not known, this trivial bound is not available. This will be an issue for designing a minimax regularization function when $\norm{t^*}_1$ is unknown and small (actually smaller than $\sigma\sqrt{\log(ed)/N}$). Somehow the ``signal-to-noise ratio" is too small for the models $\rho B_1^d$ with small $\rho$'s. Therefore, the trivial upper bound $\rho^2$ is optimal when $\rho$ is known but in the other case we will have to pay the price due to the noise and there will be no way to achieve the trivial optimal $\rho^2$ bound for small $\rho$'s (except for the trivial estimator $\hat t_0=0$, see the discussion after Proposition~\ref{prop:adapt_small_ball}). That is the reason why we will not be able to construct a minimax regularization function over the entire space $\R^d$ but only for $t^*$ such that $\norm{t^*}_1\gtrsim \sigma/\sqrt{\log(ed)/N}$. We will also show that such a construction of an optimal regularization function over the entire space $\R^d$ is actually not possible at all later in Proposition~\ref{prop:adapt_small_ball}.

Finally let us turn to the construction of a minimax regularization function for the $\ell_1^d$-norm. To that end we will use the function $\rho\geq0\mapsto r^2(\rho) = \max\big(r_Q^2(\rho), r_M^2(\rho)\big)$ as a sharp way to measure the complexity of the model $\rho B_1^d$. The main result of this article is that this function is a minimax regularization function as introduced in Definition~\ref{def:optimality}.

\begin{Theorem}\label{theo:main}
  There are absolute constants $\eta,Q,\zeta, \zeta^\prime, \Delta_0, c_0$ such that the following holds.   When $\zeta^\prime d\geq N$ or  $\zeta d\leq N$,  a minimax regularization function for the $\ell_1^d$-norm over $\R^d \backslash (\Delta_0 \sigma\sqrt{\log(ed)/N}B_1^d)$ for the confidence parameter $\delta_N=1/4$  is given by the following function: for all $\rho>0$, 
\begin{equation*}
\Psi(\rho) = c_0 r^2(\rho) 
\end{equation*}   where $r(\rho) = \max\big(r_Q(\rho), r_M(\rho)\big)$ for
\begin{equation*}
r_Q(\rho) = \inf\Big(r>0:\ell^*(\rho B_1^d\cap r B_2^{d})= Q r \sqrt{N}\Big)
\end{equation*}and 
\begin{equation*}
r_M(\rho) = \inf\Big(r>0:\sigma \ell^*( \rho B_1^d\cap r B_2^{d})= \eta r^2 \sqrt{N}\Big)
\end{equation*} denoting by $\ell^*( \rho B_1^d\cap r B_2^{d})$ the Gaussian mean width of the localized sets  $\rho B_1^d\cap r B_2^{d}$. In that case, the rate achieved by the RERM $\hat t_{\Psi}$ is the minimax rate $r^2(\norm{t^*}_1)$ when $\norm{t^*}_1\geq \Delta_0\sigma \sqrt{\log(ed)/N}$.
 \end{Theorem}

\begin{figure}[h]
\begin{center}
\begin{tabular}{cc}
   \includegraphics[width=0.5\linewidth]{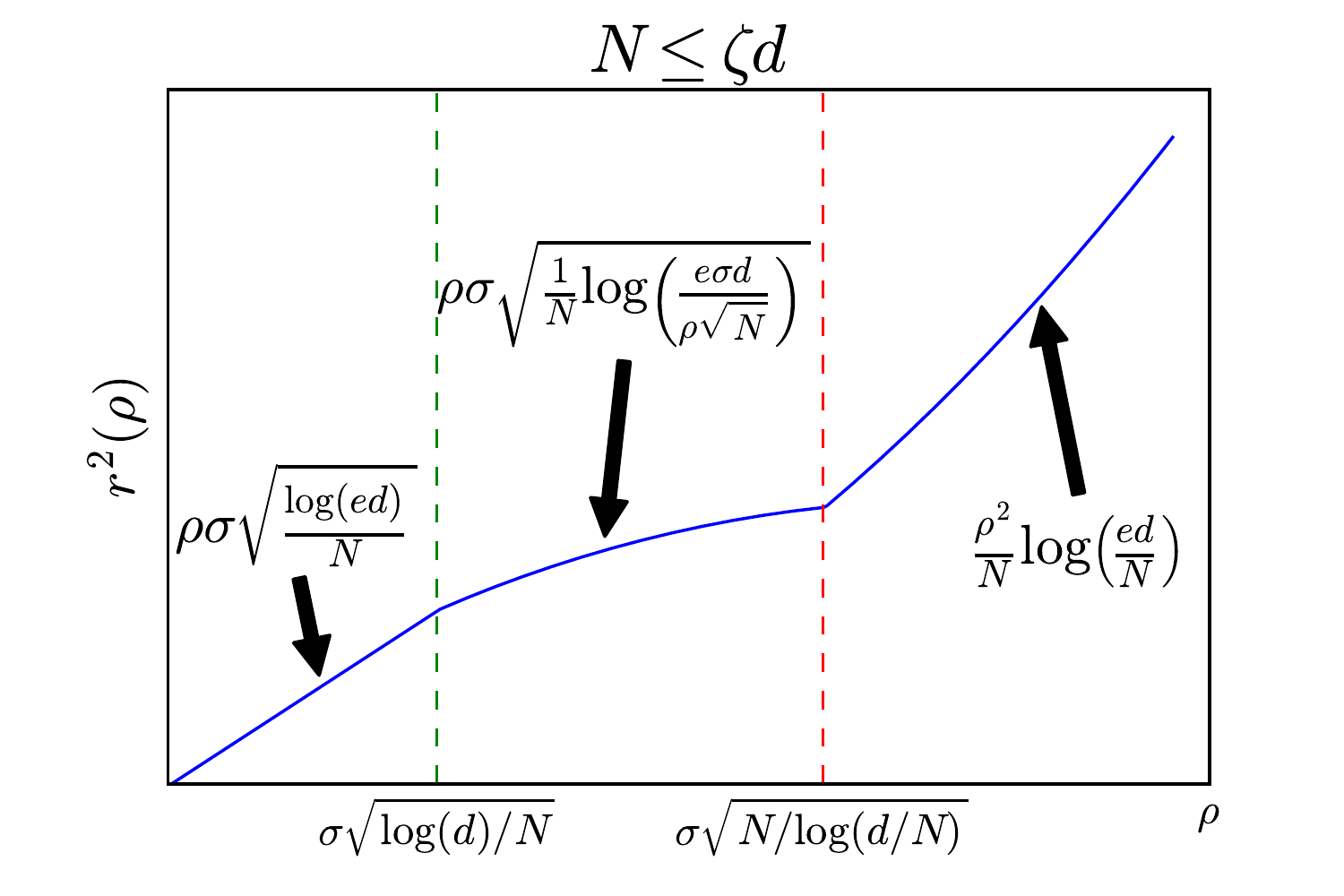} 
   \includegraphics[width=0.5\linewidth]{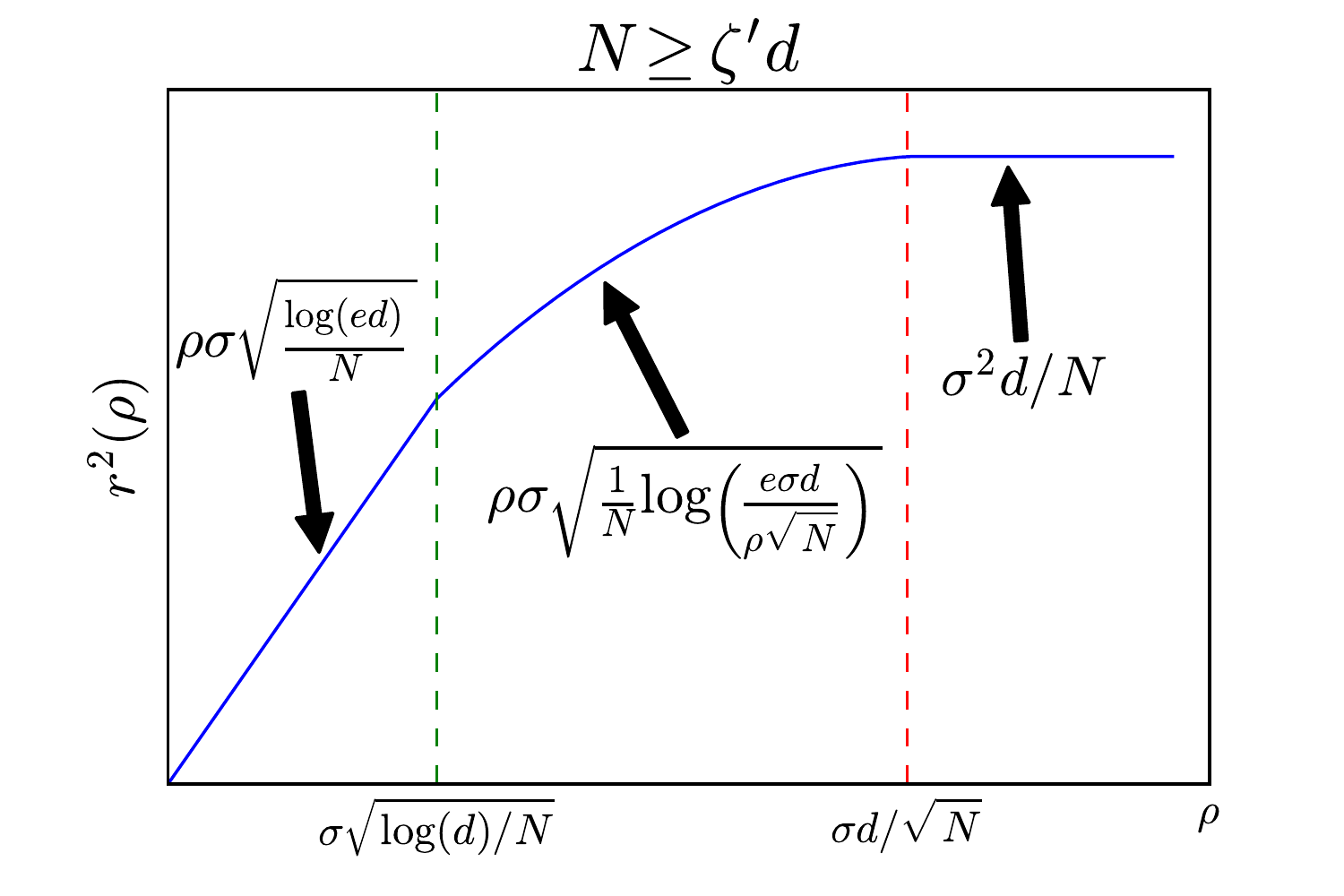} 
\end{tabular}
\caption{Shape of the graph of the minimax regularization function  $\rho\to r^2(\rho)$ of the $\ell_1^d$-norm for the cases $N\leq \zeta d$ (left) and $N\geq \zeta^\prime d$ (right)}\label{fig:mini_reg_func}
\end{center}
\end{figure}

 The shape of the minimax regularization function $\rho\to\Psi(\rho) = c_0 r^2(\rho)$ is given in Figure~\ref{fig:mini_reg_func} in the two cases $N\leq \zeta d$ (``high-dimensional statistics'') and $N\geq \zeta^\prime d$ (``classical or low-dimensional statistics''). 

The only difference between the two cases ($\zeta^\prime d\geq N$ or  $\zeta d\leq N$) appears for large radii $\rho$. The reason for that lies in the statistical complexity coming (or not) from the distortion of the operator $\bX$. In the low-dimensional case, $\bX$ is such that (with high probability), $\norm{\bX t}_2 \sim \norm{t}_2$  for all $t\in\R^d$. There is no distortion coming from $\bX$. Somehow observing $\bX t^*$ is the same as observing $t^*$ itself, one just has to invert $\bX$ -- this can be done because $\bX$ acts like an isomorphy on the entire space $\R^d$. Therefore, there is no statistical complexity coming from $\bX$ and so its associated complexity parameter $r_Q(\cdot)$ does not show up in the final complexity parameter $r(\cdot)=\max(r_M(\cdot), r_Q(\cdot))$. We therefore end up with $r(\cdot) = r_M(\cdot)$ in the low-dimensional case. In particular, for large radii $\rho$ (for which, one has $\rho B_1^d\cap r_M(\rho) B_2^d = r_M(\rho) B_2^d$), we pay the worst rate of convergence in $\R^d$, which is $\sigma^2 d/N$ because learning over $\rho B_1^d$ for large values of $\rho$ is as hard as learning over the entire space $\R^d$ and the price for the latter is the rate $\sigma^2 d/N$. 

 The situation is totally different in the high-dimensional setup because in that case the operator $\bX\in\R^{N\times d}$ has a none trivial kernel; therefore, observing $\bX t^*$ is totally different from observing $t^*$ (for instance, imaging that $t^*\in \ker \bX$). This adds to the statistical complexity of the problem of estimating $t^*$. In this regime, both the noise and the distortion effect of $\bX$ appear in the statistical complexity of the estimation problem; this means that both complexity parameter $r_Q(\cdot)$ and $r_M(\cdot)$ appear in the total complexity parameter $r(\cdot)$ and therefore in the ultimately designed minimax regularization function. For small values of $\rho$, the effect of the noise is predominant but for large values of $\rho$ this is the effect of $\bX$ which is the main responsible of the statistical complexity. In particular, the $\ell_2^d$-diameter of $\ker \bX \cap \norm{t^*}_1 B_1^d$ is important because there is of course no way to distinguish $t^*$ from $t^*+h$ for all $h\in\R^d$ such that $\bX t^* = \bX(t^*+h)$ that is for all $h\in \ker \bX$ such that $t^*+h\in \norm{t^*}_1 B_1^d$. Hence, estimating $t^*$ is at least as hard as estimating any point in $(t^*+\ker \bX)\cap \norm{t^*} B_1^d$ and therefore, no estimator $\tilde t$ can estimate $t^*$ at a rate better than ${\rm diam}\left(\ker \bX \cap \norm{t^*}_1 B_1^d, \ell_2^d\right)^2$. The latter quantity is itself lower bounded by the Gelfand's $N$-width of $\norm{t^*} B_1^d$ defined as
 \begin{equation}\label{eq:gelfand}
  c_N(\norm{t^*}_1 B_1^d):= \inf\left\{{\rm diam}(\ker \Gamma \cap \norm{t^*}_1 B_1^d) : \Gamma \in\R^{N\times d}\right\} \sim \norm{t^*} \min\left\{1, \sqrt{\frac{\log(ed/N)}{N}}\right\} 
  \end{equation} the latter result is due to Garanaev and Gluskin \cite{MR759962}. It appears that the Gelfand's $N$-width of $\norm{t^*} B_1^d$ are achieved (up to absolute constants and with high probability) by the kernel of standard $N\times d$ Gaussian matrices, which is exactly the case of the design matrix $\bX$. Therefore, with high probability,
  \begin{equation}\label{eq:third_regime}
  {\rm diam}\left(\ker \bX \cap \norm{t^*}_1 B_1^d, \ell_2^d\right)^2\sim c_N^2(\norm{t^*}_1 B_1^d) \sim \frac{\norm{t^*}_1}{N}\log\left(\frac{ed}{N}\right) \sim r_Q^2(\norm{t^*}_1)
  \end{equation} This is exactly the price we pay in $r_Q(\rho)$ when $\rho\geq \sigma \sqrt{N/\log(ed/N)}$. That is the reason why we take the regularization function $\Psi(\rho)$ of the order of the Gelfand's $N$-width of $\rho B_1^d$ (to the square) for large radii $\rho$: it is the right concept of statistical complexity that shows up in this part of the space $\R^d$, where the statistical complexity coming from the distortion of $\bX$ becomes more important than the one due to the noise.

\begin{Remark}[Regularization function for the LASSO]
The LASSO is the RERM procedure obtained for a linear regularization function $\Psi(\rho) = \sigma \rho \sqrt{\log d / N}$ which is obtained by using a trivial upper bound on the complexity of the model $\rho B_1^d$ in \eqref{eq:minimax_lasso}: 
\begin{equation}\label{eq:lasso_maj}
{\rm comp} \left(\rho B_1^d \right) = r_M^{2}(\rho) \leq \frac{\sigma \ell^*\left(\rho B_1^d\right)}{\sqrt{N}}  = \sigma \rho \sqrt{\frac{\log(ed)}{N}}.
\end{equation}
 This complexity is obtained by simply removing the localization (i.e. the intersection with $rB_2^d$) in the multiplier process when computing $r_M(\cdot)$, and does not take $r_Q(\cdot)$ into account. This means that the distortion of the operator $\bX$ is supposed to have no effect on the statistical complexity of the problem. This is why estimation results for the LASSO deal only with the reconstruction of vectors which are sparse or almost sparse, i.e. for vectors belonging to the cone appearing in the RE or CC conditions, cf. \cite{MR2533469}. Over this cone, the quadratic process behaves nicely (that is, the isomorphic property from Proposition~\ref{prop:event_Omega^*} holds on this cone) or in other words, the operator $\bX$ is well-conditioned on the set of vectors we want to reconstruct, so that there is no statistical complexity coming from the distortion of this operator. So, as long as estimation of sparse or approximately sparse vectors is concerned, there is no need for the complexity function $r_Q(\cdot)$. That is why the regularization function used for the LASSO take into account only the fixed point $r_M(\cdot)$ associated to the statistical complexity due to the noise and not the one from the inverse problem. On the contrary, by taking $r_M(\cdot)$ and $r_Q(\cdot)$ into account, our regularization function allows us to deal with the full space $\mathbb{R}^d$ (except for a small $\ell_1^d$-ball centered in $0$, cf. Proposition~\ref{prop:adapt_small_ball})  and not only a cone.
 
   Moreover, as said before, the way the regularization function is designed in \eqref{eq:lasso_maj} is sub-optimal because it uses a trivial upper bound on $r_M(\cdot)$ instead of using the exact formulation of $r_M(\cdot)$ as in \eqref{eq:rM}. Contrary to the LASSO, this latter exact formulation takes into account, thanks to the localization, the fact that the regularization is not needed on the whole space --in some areas the random processes behave nicely whatever. The suboptimal approach for the LASSO is likely to be responsible for a loss in the rate of convergence achieved by the LASSO, which is $ s \sigma^2\log(ed)/N$ whereas the minimax rate is $s\sigma^2\log(ed/s)/N$ (cf. \cite{BLT}). This is not a big loss, especially when $d>>s$, but from a purely theoretical point of view the right way to regularize for the reconstruction of sparse vectors should be using  $r_M^2(\norm{t}_1)$ instead of $\sigma \norm{t}_1 \sqrt{\log(ed)/N}$ as it is the case for the LASSO. However, the resulting regularization function would be concave (cf. the right-hand side plot in Figure~\ref{fig:mini_reg_func}). Therefore, the small price paid from a theoretical point of view by using the trivial upper bound in \eqref{eq:lasso_maj} seems to be worth the computational gain obtained by using a convex regularization as does the LASSO.
\end{Remark}

 Let us now turn to the adaptation problem in the ball $\rho B_1^d$ for $\rho\sim \sigma \sqrt{\log(ed)/N}$. We want to answer the following question: is it possible to construct a regularization function $\Psi(\cdot)$ so that the associated regularized procedure $\hat t_\Psi$ is adaptive on the entire space $\R^d$? Or (even stronger) is there any statistic that can be adaptive (in the sense that it achieves the rate of the ERM on $\norm{t^*}_1B_1^d$ without knowing $\norm{t^*}_1$ beforehand) on the entire space $\R^d$ (this statistics may not be a regularized procedure)? It appears that the answer to this question is negative, which we prove in the following proposition. 

 However, one needs to be cautious with the next statement because there is a trivial estimator $\hat t_0 = 0$ such that for every $t^*\in\R^d$, with probability $1$, $\norm{\hat t_0 - t^*}_2^2 = \norm{t^*}_2^2$ and therefore $\hat t_0$ is adaptive on $\rho B_1^d$ as long as the minimax rate over $\rho B_1^d$ is $\rho^2$,  which is the case for any $\rho\lesssim \sigma \sqrt{\log(ed)/N}$. Therefore, there exists a procedure adaptive on $\rho B_1^d$ when $\rho\sim \sigma \sqrt{\log(ed)/N}$. Moreover, according to Theorem~\ref{theo:main}, there exists a procedure adaptive on  $\R^d\backslash  \rho B_1^d$. But the question concerns the adaptation on the \textit{entire space $\R^d$ at the same time}.  

 The following statement shows that if $\hat t$ is a procedure adaptive on $\R^d\backslash \rho B_1^d$ then it cannot be adaptive on $\rho B_1^d$ for $\rho \sim \sigma \sqrt{\log(ed)/N}$. Moreover, it also proves that adaptation on the entire space $\R^d$ is not possible and that Theorem~\ref{theo:main} is optimal given that the range of radii $[\Delta_0\sigma \sqrt{\log(ed)/N}, +\infty)$ on which it is adaptive cannot be inflated (up to absolute constants). Before turning to the statement let us denote by $\bP_{t^*}$ the probability distribution of a $N$-sample $(X_i, Y_i)_{i=1}^N$ of i.i.d. copies of $(X, Y)$ when $(X, Y)$ is distributed according to \eqref{eq:model}.

 \begin{Proposition}\label{prop:adapt_small_ball}
 Assume that $2d\geq \exp(544/225)$ and that there exists an absolute constant $\chi_1$ such that the following holds. Let $\rho  \leq 2 \sigma\sqrt{(\log(2d))/(96N)}$ be such that $16 \chi_1 r^2(\rho)\leq \rho^2$ and denote by $(e_j)_{j=1}^d$ the canonical basis of $\R^d$. Assume that $\hat t$ is an estimator such that for every $t^*\in\{\pm \rho e_1, \ldots, \pm \rho e_d\}$, 
\begin{equation*}
\bP_{t^*}\left[\norm{\hat t - t^*}_2^2\leq \chi_1 r^2(\rho)\right]\geq \frac{3}{4}.
\end{equation*}Then, for every $t^*\in  (\rho/2) B_1^d$, 
\begin{equation}\label{eq:result}
\bP_{t^*}\left[\norm{\hat t - t^*}_2^2\geq \rho^2/16\right]\geq \frac{1}{2}.
\end{equation}
 \end{Proposition}

The proof of Proposition~\ref{prop:adapt_small_ball} is given in Section~\ref{sec:technical_material}. Note that the only property of the design $X$ used to prove Proposition~\ref{prop:adapt_small_ball} is isotropicity. Since isotropicity does not tell much on the distortion properties of the design matrix $\bX$, it means that Proposition~\ref{prop:adapt_small_ball} is only based on the statistical complexity coming from the noise. This is not a surprise given that Proposition~\ref{prop:adapt_small_ball} is a result for very small radii less than $\sim \sigma \sqrt{\log(ed)/N}$. At that scale, even if $\ker \bX$ is in the worst possible position, i.e.  ${\rm diam}(\ker \bX \cap \rho B_1^2, \ell^d_2)^2 = {\rm diam}( \rho B_1^2, \ell^d_2)^2 = \rho^2$, we still have $\rho^2\lesssim r_M^2(\rho)$. Hence, the distortion of $\bX$ does not play any role at this very small scale and therefore that is not a surprise that Proposition~\ref{prop:adapt_small_ball} is true for any isotropic design $X$.  

Finally, let us rephrase Proposition~\ref{prop:adapt_small_ball} in other words. Proposition~\ref{prop:adapt_small_ball} shows that if a procedure can learn all vectors in $\{\pm \rho e_1, \ldots, \pm \rho e_d\}$ at the minimax rate $r^2(\rho)$ then this estimator cannot learn any $t^*\in (\rho/2)B_1^d$ at the optimal minimax rate $\rho^2$ for confidence $1/4$. For instance, given that the result \eqref{eq:result} holds for any $t^*\in (\rho/2)B_1^d$, in particular, for $t^*=0$, it tells that $\hat t$ cannot estimate $t^*=0$ at a rate better than $\rho^2\sim \sigma^2 \log(ed)/N$ whereas the minimax rate over $\rho^* B_1^d$ for $\rho^*=0$ is obviously $0$. Finally, note that the condition $16 \chi_1 r^2(\rho)\leq \rho^2$ implies that $\rho\gtrsim \sigma \sqrt{\log(ed)/N}$ so that the phase transition radius above which adaptation is possible but not below is of the order of $\sigma \sqrt{\log(ed)/N}$ which is the radius we have found in Theorem~\ref{theo:main}.


\section{Proof of Theorem~\ref{theo:main}}
 \label{sec:proof_of_theorem_ref_theo_main}
 
Most of the proof consists in showing that with high probability, $\hat t$ belongs to  $t^*+\rho^* B_1^d$ where 
\begin{equation}\label{eq:rho_star}
\rho^*=\max\left(10,8\frac{C_M^{(2)})^2}{(C_M^{(1)})^2\eta}+1\right) \norm{t^*}_1.
\end{equation}
(once this goal is achieved, it is straightforward to show (again with high probability) that  
$\norm{\hat{t}-t^*}_2^2$ is less than the minimax rate of convergence over $\rho^* B_1^d$).

To do so, we will prove that with high probability, any $t$ outside $t^*+\rho^* B_1^d$ satisfies $P_N\cL_t^\Psi>0$ (whereas $P_N\cL_{\hat t}^\Psi \leq 0$). We partition $\R^d\backslash (t^*+\rho^* B_1^d)$ into shelves of the form $t^*+(2^{j+1}\rho^* B_1^d \backslash (2^j\rho^*) B_1^d)$, in which the regularization function remains mostly constant.  We only need to study the smallest shells, i.e. for $k=1, \ldots, K_0$ for some well-chosen $K_0$ ($K_0$ is the smallest integer so that $2^{K_0-1}\rho^* B_1^d\cap r(2^{K_0-1}\rho^*) B_2^d = r(2^{K_0-1}\rho^*) B_2^d$), the part of $\R^d$ for which $\norm{t}_1\geq 2^{K_0}\rho^*$ will be treated by an homogeneity argument. 

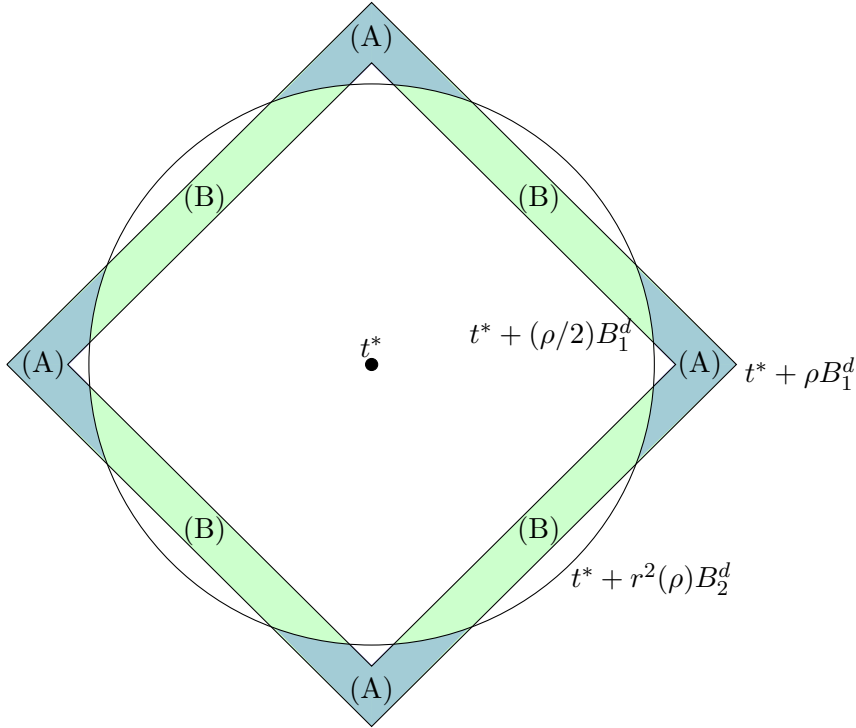
\begin{figure}[!h]
 \centering
 \begin{tikzpicture}[scale=0.4]

 \fill[color=green, opacity=0.2,even odd rule] (-12,0) -- (0,12) -- (12,0) -- (0,-12) -- (-12,0) -- (-10,0) -- (0,10) -- (10,0) -- (0,-10) -- (-10,0) -- cycle ;
   \fill[color=blue,opacity=0.2] (12,0) -- (9.3,0) --(9.3,0) arc (0:20:9.3) -- cycle ;
     \fill[color=blue,opacity=0.2] (12,0) -- (9.3,0) --(9.3,0) arc (0:-20:9.3) -- cycle ;
        \fill[color=blue,opacity=0.2] (-12,0) -- (-9.3,0) --(-9.3,0) arc (0:20:-9.3) -- cycle ;
     \fill[color=blue,opacity=0.2] (-12,0) -- (-9.3,0) --(-9.3,0) arc (0:-20:-9.3) -- cycle ;
        \fill[color=blue,opacity=0.2] (0,12) -- (0,9.3) --(0,9.3) arc (90:110:9.3) -- cycle ;
     \fill[color=blue,opacity=0.2] (0,12) -- (0,9.3) --(0,9.3) arc (90:70:9.3) -- cycle ;     
        \fill[color=blue,opacity=0.2] (0,-12) -- (0,-9.3) --(0,-9.3) arc (90:110:-9.3) -- cycle ;
     \fill[color=blue,opacity=0.2] (0,-12) -- (0,-9.3) --(0,-9.3) arc (90:70:-9.3) -- cycle ;    
      \fill[color=white, opacity=1.0] (0,-9.95) -- (-9.95,0) -- (0,9.95) -- (9.95,0) -- cycle ;
   \filldraw (0,0) circle (0.2cm);
 \draw (0,0.6) node{$t^*$};
 \draw  (-10,0) -- (0,10) -- (10,0) -- (0,-10) -- (-10,0);
 \draw (0,0) circle (9.3cm);
 \draw (-12,0) -- (0,12) -- (12,0) -- (0,-12) -- (-12,0);    
      
  \draw (10.8,0) node {(A)};  
   \draw (-10.8,0) node {(A)};
    \draw (0,10.8) node {(A)};
     \draw (0,-10.8) node {(A)}; 
      \draw (5.5,-5.5) node {(B)}; 
      \draw (5.5,5.5) node {(B)}; 
      \draw (-5.5,-5.5) node {(B)}; 
      \draw (-5.5,5.5) node {(B)}; 
       \draw (14.1,-0.3) node {$t^*+\rho B_1^d$};
       \draw (5.9,1) node {$t^*+(\rho/2) B_1^d$};   
        \draw (9.2,-7.1) node {$t^*+r^2(\rho) B_2^d$};

 \end{tikzpicture}
 \caption{$P_N {\cal M}_{t-t^*}$ is dominated by $P_N {\cal Q}_{t-t^*}$ in (A) and by  ${\cal R}_{t,t^*}$ in (B). The regularization function $\Psi(\cdot)$ is designed in order to dominate $P_N {\cal M}_{t-t^*}$ in (B).}
 \label{fig:intuition}
 \end{figure}

On each of the smallest shelves, the argument is roughly and heuristically the following: we place ourselves on a high probability event on which random processes ($P_N {\cal M},P_N {\cal Q}$, their supremum, their infimum,...) ``behave  nicely" (i.e. they both scale like $\norm{t-t^*}_2^2$). Then, as $P_N {\cal M}_{t-t^*}$ is the only possibly negative term, it suffices to identify zones where $P_N {\cal Q}_{t-t^*}>P_N {\cal M}_{t-t^*}$ (then directly $P_N\cL_t^\Psi>0$), and compensate $|P_N {\cal M}_{t-t^*}|$ on the other part by using a penalty that is close to the supremum (on this other part) of $P_N {\cal M}_{t-t^*}$. As $P_N {\cal Q}_{t-t^*}$ grows quicker than $P_N {\cal M}_{t-t^*}$ with respect to $\norm{t-t^*}_2$, the big zone where $P_N {\cal Q}_{t-t^*}>P_N {\cal M}_{t-t^*}$ will be the exterior of $t^*+rB_2^d$ for an adequate $r$ (cf. Figure~\ref{fig:intuition}).  This $r$ must be such that any $t$ in the exterior of this ball satisfies $P_N {\cal Q}_{t-t^*}\gtrsim \norm{t-t^*}_2^2 \gtrsim P_N {\cal M}_{t-t^*}$, and we will see that the first inequality amounts to $r\geq r_Q(\rho)$, and the second one to $r\geq r_M(\rho)$. Next, the supremum of $P_N {\cal M}_{t-t^*}$ on $\rho B_1^d\cap r(\rho)B_2^d$ is less than $r_M(\rho)r(\rho)\leq r^2(\rho)$ for $r(\rho)=\max(r_M(\rho),r_Q(\rho))$. We therefore set the regularization function $\Psi(\rho)$ at level $\rho$ to be proportional to the quantity  $r^2(\rho)$ because it is this quantity measuring the amplitude of the oscillation of the multiplier process in $\rho B_1^d \cap r(\rho)B_2^d$.
 
As for its presentation, the proof of Theorem~\ref{theo:main} is divided into two parts. The first part (Section~\ref{sec:probabilistic}) defines the event on which the two processes ``behave nicely" and computes a lower bound on its probability. In the second part (Section~\ref{sec:deterministic}) we will place ourselves on this event and carry out the deterministic geometric part of the argument.

\subsection{Probabilistic control of the processes}
\label{sec:probabilistic}
Instead of controlling the two processes on shelves, we will control them on the full $\ell_1$ balls, because it does not change the complexity, up to constants, and the very last step of the proof requires a control on the two processes on the full $\ell_1$ ball $\rho^* B_1^d$. 


\subsubsection{Control of the quadratic process}
This first section provides the classical analysis of the quadratic process based upon its isomorphic properties on the set of ``almost sparse vectors''. Such a property holds in the optimal regime of observation (or the optimal size of the cone of ``almost sparse vectors''), only in the sub-Gaussian case.  It  is the case we are considering here since we assumed that the design is a standard Gaussian random variable. This analysis borrows some ideas from the ``isomorphic method'' from \cite{MR2240689} or the Restricted Isometry Property from \cite{MR2236170} in the sub-Gaussian case. For the sake of completeness we recall here the argument from \cite{LM13}.

\begin{Proposition}\label{prop:event_Omega^*}There are absolute constants $ C_1$ and $C_1^\prime$ such that the following holds. Let $X_1,\ldots, X_N$ be $N$ i.i.d. standard Gaussian vectors in $\R^d$. Denote by $\Omega^*$ the event on which: for every $\rho \geq \rho^*$ and all  $t\in t^*+\rho B_1^d$,
\begin{equation}\label{eq:isomorphy}
\mbox{ if } \norm{t-t^*}_2\geq r_Q(\rho) \mbox{ then } \frac{1}{2}\norm{t-t^*}_2^2 \leq \frac{1}{N}\sum_{i=1}^N \inr{X_i, t-t^*}^2 \leq \frac{3}{2}\norm{t-t^*}_2^2.
\end{equation}Then, one has  $\bP[\Omega^*]\geq1-2\exp(-C_1 Q^2N)$ as long as $Q\leq C_1^\prime$. 
\end{Proposition}
\beginproof
First note that for all $\rho>0$,  $r_Q(\rho)=\rho r_Q(1)$. Indeed, we have
$$ \ell^*(\rho B_1^d \cap rB_2^d)=\ell^*\left(\rho(B_1^d \cap (r/\rho)B_2^d)\right)=\rho \ell^*\left( B_1^d \cap (r/\rho)B_2^d\right)$$
and so  
\begin{align}\label{eq:homo_r_Q}
 \nonumber r_Q(\rho)&=\inf\{r>0: \ell^*(\rho B_1^d \cap rB_2^d)  = Q \sqrt{N} r\}  \\
\nonumber &=\inf\{r>0: \ell^*(  B_1^d \cap (r/\rho)B_2^d)   =  Q \sqrt{N} (r/\rho)\} \\
&=\rho \  \inf\{r>0: l^*( B_1^d \cap rB_2^d)   =  Q \sqrt{N} r\}  = \rho r_Q(1).
\end{align}

For all $\rho>0$, define the event $\Omega(\rho)$ on which one has for all $t\in t^*+\rho B_1^d$, 
\begin{equation*}
\mbox{ if } \norm{t-t^*}_2\geq r_Q(\rho) \mbox{ then } \frac{1}{2}\norm{t-t^*}_2^2 \leq \frac{1}{N}\sum_{i=1}^N \inr{X_i, t-t^*}^2 \leq \frac{3}{2}\norm{t-t^*}_2^2.
\end{equation*}
Let us show that if $\Omega(\rho^*)$ holds then for any $\rho\geq \rho^*$, $\Omega(\rho)$ holds as well. Suppose that $\Omega(\rho^*)$ holds. Consider $t \in t^*+\rho B_1^d\textup{ such that }\norm{t-t^*}_2>r_Q(\rho)$ and define 
$$t^\prime:=t^*+(\rho^*/\rho)(t-t^*) \in t^*+\rho^*B_1^d.$$ It follows from \eqref{eq:homo_r_Q} that  $r_Q(\rho^*)=(\rho^*/\rho)r_Q(\rho)$.  Thus $\norm{t^\prime-t^*}_2=(\rho^*/\rho)\norm{t-t^*}_2>(\rho^*/\rho) r_Q(\rho)=r_Q(\rho^*)$, and since $\Omega(\rho^*)$ holds, it follows that $\norm{t^\prime-t^*}_2^2/2 \leq \sum_{i=1}^N \inr{X_i, t^\prime-t^*}^2 /N\leq 3\norm{t^\prime-t^*}_2^2/2$. This implies that $\norm{t-t^*}_2^2/2 \leq \sum_{i=1}^N \inr{X_i, t-t^*}^2 /N\leq 3\norm{t-t^*}_2^2/2$ so $\Omega(\rho)$ holds.

As a conclusion, $\Omega^*=\Omega(\rho^*)$ and we can now lower bound the probability that this event holds.

Let us consider the class of linear functions
 $$F =\left\{ \inr{\cdot,t-t^*}, t \in t^*+ \rho^* B_1^d \cap r_Q(\rho^*)S_2^{d-1}\right\} = \left\{ \inr{\cdot,t}, t \in \rho^* B_1^d \cap r_Q(\rho^*)S_2^{d-1}\right\}.$$ 

We assume that $F$ is non empty (if $F=\emptyset$ then the theorem is trivially satisfied). It follows from Theorem~1.12 in \cite{shahar_multi_pro} that for any $x>0$, with probability at least $1-2 \exp\big(-C_1 \min(x^2, x\sqrt{N})\big),$
\begin{equation*}
\sup_{f\in F}\left|\frac{1}{N}\sum_{i=1}^N f^2(X_i) - \E f^2(X)\right|\leq C_2\left(\frac{\Delta\gamma}{\sqrt{N}}+\frac{\gamma^2}{N} + \frac{x \Delta^2}{\sqrt{N}}\right)
\end{equation*}where $\Delta$ is the diameter in $\psi_2$ of $F$ and $\gamma$ is Talagrand's $\gamma_2$ functional of $F$ w.r.t. $\psi_2$. Note that since $X$ is a standard Gaussian variable in $\R^d$, for any $t\in\R^d$,  $\norm{\inr{X, t}}_{\psi_2} = C_3 \norm{t}_{2}$ for some absolute constant $C_3$. It follows that
\begin{equation*}
\Delta=2 \sup_{t\in \rho^* B_1^d \cap r_Q(\rho^*)S_2^{d-1}}C_3\norm{t}_{2} = 2C_3r_Q(\rho^*) \mbox{ and } \gamma = \gamma_2(\rho^* B_1^d \cap r_Q(\rho^*)S_2^{d-1}, \ell_2^d). 
\end{equation*} Moreover, it follows from the Majorizing measure theorem (cf. Chapter~1 in \cite{MR3184689}) that 
\begin{equation*}
\gamma_2(\rho^* B_1^d \cap r_Q(\rho^*)S_2^{d-1}, \ell_2^d)\leq C_4 \ell^*(\rho^* B_1^d \cap r_Q(\rho^*)S_2^{d-1})
\end{equation*} 
Since $F$ is non-empty, by Lemma~\ref{prop: equivalent localizations} the right-hand side is equal to $C_4 \ell^*(\rho^* B_1^d \cap r_Q(\rho^*)B_2^{d-1})$  and so by definition of $r_Q(\rho^*)$, one has $\gamma\leq C_4 Q r_Q(\rho^*)\sqrt{N}$.

Since $X$ is isotropic (i.e. for any $t\in\R^d$, $\E\inr{X, t}^2 = \norm{t}_2^2$), we obtain for $x=Q\sqrt{N}$ that, with probability greater than $1 - 2\exp(-C_1 \min(Q,Q^2)N)$, for any $t \in t^*+ \rho^* B_1^d \cap r_Q(\rho^*)S_2^{d-1}$, 
\begin{equation*}
\left| \frac{1}{N}\sum_{i=1}^N \inr{X_i, t-t^*}^2 -\norm{t-t^*}_2^2 \right| \leq (2C_2C_3C_4Q+C_2C_4^2Q^2+4C_2QC_3^2)r_Q^2(\rho^*)
\end{equation*} 

So, as long as long as: $Q\leq C_1^\prime:= \min\{1,(12C_2C_3C_4)^{-1},(\sqrt{6C_2}C_4)^{-1},(24C_2C_3^2)^{-1}\}$, one has, with probability greater than $1 - 2\exp(-C_1 Q^2N)$,  for all $t \in t^*+ \rho^* B_1^d \cap r_Q(\rho^*)S_2^{d-1}$,
\begin{equation}\label{eq:isom_prop_sphere}
 \left| \frac{1}{N}\sum_{i=1}^N \inr{X_i, t-t^*}^2 -\norm{t-t^*}_2^2 \right|\leq \frac{r_Q^2(\rho^*)}{2}=\frac{1}{2}\norm{t-t^*}_2^2.
\end{equation} In other words, the quadratic process satisfies an isomorphic property on the set $t^*+ (\rho^* B_1^d \cap r_Q(\rho^*) S_2^{d-1})$. Now, it remains to extend this result to the set of vectors $t\in t^*+\rho^* B_1^d$ such that $\norm{t-t^*}_2\geq r_Q(\rho^*)$. Let $t$ be such a vector and define $t^\prime:= t^* + (r_Q(\rho^*)/\norm{ t-t^* }) (t-t^*)$. Since  $t^\prime\in t^*+ (\rho^* B_1^d \cap r_Q(\rho^*) S_2^{d-1})$, it satisfies the isomorphic property from \eqref{eq:isom_prop_sphere} and so
$$ \frac{1}{2}\norm{t-t^*}_2^2 \leq  \frac{1}{N}\sum_{i=1}^N\inr{X_i, t-t^*}^2 \leq \frac{3}{2}\norm{t-t^*}_2^2$$
which corresponds exactly to the event $\Omega(\rho^*)$. Therefore, $\bP[\Omega^*] = \bP[\Omega(\rho^*)]\geq 1-2\exp(-C_1Q^2 N)$.
\endproof

\subsubsection{Control of the multiplier process}
In this section, we provide a control of the multiplier process on several shelves of the space $\R^d$.  

Define $r_0$ as the non-null solution to $\sigma \ell^*(rB_2^{d})=\eta r^2 \sqrt{N}$ -- i.e. $r_0 = \sigma \ell^*(B_2^d)/(\eta \sqrt{N})= (\sigma/\eta)\sqrt{d/N}$. Let  $\rho_0$  be the smallest $\rho$ such that $\rho B_1^d$ contains $r_0B_2^{d}$ --i.e., $\rho_0=r_0\sqrt{d}$. We can see that $\rho_0$ is such that $r_M(\rho) = r_M(\rho_0)$ for all $\rho\geq \rho_0$. Indeed, one first sees that $r_M(\rho_0)=r_0$, since $\sigma \ell^*(r_0 B_2^d\cap \rho_0 B_1^d)= \sigma \ell^*(r_0 B_2^d) = \eta r_0^2 \sqrt{N}$, and $\sigma \ell^*(r B_2^d\cap \rho_0 B_1^d) = \sigma r \ell^*(B_2^d)>\eta r^2 \sqrt{N}$ for all $r< r_0$  and for $r>r_0$, $\sigma \ell^*(rB_2^{d} \cap \rho_0 B_1^d)\leq \sigma r \ell^*(B_2^d)\leq \eta r^2 \sqrt{N}$. This last argument also holds for $\rho \geq \rho_0$. In the latter case, if $r>r_0$ then
$$\sigma \ell^*(rB_2^{d} \cap \rho B_1^d)\leq \frac{r}{r_0}\sigma \ell^*(r_0B_2^{d} \cap \rho B_1^d)=\frac{r}{r_0}\sigma \ell^*(r_0B_2^{d} \cap \rho_0 B_1^d)=\frac{r}{r_0} \eta r_0^2\sqrt{N}  < \eta r^2\sqrt{N}$$
which means that $r_M(\rho) \leq r_0$. And as $\rho \geq \rho_0$, $r_M(\rho) \geq r_M(\rho_0)=r_0$. Therefore, for $\rho \geq \rho_0$, $r_M(\rho)$ is constant, equal to $r_0$.
And on $[0,\rho_0]$, $r_M$ is  non-decreasing: let $\rho^\prime \leq \rho^{\prime \prime} \leq \rho_0$, then $\sigma \ell^*(r_M(\rho^\prime) B_2^{d} \cap \rho^{\prime \prime} B_1^d) \geq \sigma \ell^*(r_M(\rho^\prime) B_2^{d} \cap \rho^\prime B_1^d)=\eta r_M(\rho^\prime)^2\sqrt{N}$ so $ r_M(\rho^{\prime \prime}) \geq r_M(\rho^\prime)$.

 We denote $K_0=\min \{k \in \mathbb{N}: 2^k \rho^* \geq 2\rho_0\}$: we will see later that $K_0$ is defined that way to be the number of the first ``shell" such that $r_M(2^{K_0-1}\rho^*)B_2^d\subset 2^{K_0-1}\rho^* B_1^d$.

\begin{Proposition}\label{prop:event_multiplicatif}
There exists an absolute constant $C_5$ such that the following holds. Let $X_1,\ldots,X_N$ be $N$ i.i.d. standard Gaussian vectors in $\R^d$ and $\xi_1,\ldots,\xi_N$ be $N$ standard real-valued Gaussian variables independent of the $X_i$'s. For all $k=0,\ldots,K_0$, denote by $A_k$ the event on which, for every  $t\in\R^d$ such that  $\norm{t-t^*}_1\leq 2^{k}\rho^* $:
\begin{equation}\label{eq:multiplier_control}
\left|P_N {\cal M}_{t-t^*}\right| \leq \frac{1}{4} \max\left(r_M(2^{k}\rho^*)^2, \norm{ t-t^* }_2^2\right) .
\end{equation}
Then,  for $\eta= 1/(16\sqrt{2})$, one has  
$$\bP\left[\bigcap \limits_{k=0}^{K_0} A_k\right] \geq 1-2\exp \left( -C_5 N  \right) - 40 \exp \left( - \frac{C_M^{(1)}Nr_M( \rho^*)^2)}{1024 C_M^{(2)}\sigma^2}\right)$$
 when $\rho^*\geq 4096 \log(2)  \sigma /\left(C_M^{(1)} \sqrt{N}\right).$
\end{Proposition}

\beginproof
We first work conditionally to the $\xi_i, i=1, \ldots, N$. Let $\rho>0$ and define $T(\rho):=t^*+\rho B_1^d \cap r_M(\rho)B_2^{d}$.  It follows from the Gaussian concentration inequality (cf. Borell's inequality in \cite{MR1849347}) that, for all $x>0$, with probability greater than $1-2\exp(-x^2/2)$,  
$$\left|\sup_{t \in T(\rho)} \sum_{i=1}^{N} \xi_i\inr{X_i,t-t^*}-\E\sup_{t \in T(\rho)}\sum_{i=1}^{N} \xi_i \inr{X_i,t-t^*}\right|\leq x \sigma(T(\rho)) $$ 
where $$ \sigma(T(\rho))= \sup \limits_{t \in T(\rho)} \sqrt{\E  \left(\sum_{i=1}^{N} \xi_i\inr{X_i,t-t^*}\right)^2 } .$$

Conditionally to the $\xi_i$'s, the Gaussian process $\left(\sum_{i=1}^{N} \xi_i\inr{X_i,t-t^*}\right)_{t\in T(\rho)}$ has the same distribution as the Gaussian process $\left(\widehat{\sigma_N}\inr{X_1,t-t^*}\right)_{t\in T(\rho)}$ for $\widehat{\sigma_N}:=\sqrt{\sum_{i=1}^{N} \xi_i^2}$.
 This yields  
\begin{align*}\E\left[\sup \limits_{t \in T(\rho)} \sum_{i=1}^{N} \xi_i \inr{X_i,t-t^*}\right]&=\widehat{\sigma_N} \ell^*(T(\rho)-t^*) =\widehat{\sigma_N} \ell^*(\rho B_1^d \cap r_M(\rho)B_2^{d})=\widehat{\sigma_N} \frac{\eta  \sqrt{N} r_M^2(\rho)}{\sigma}
\end{align*}and $\sigma(T(\rho)) = \sup_{t \in T(\rho)} \widehat{\sigma_N} \sqrt{\E[\inr{X,t-t^*}^2]} = \widehat{\sigma_N} \sup_{t\in T(\rho)}\norm{t-t^*}_2\leq \widehat{\sigma_N} r_M(\rho)$. 

So, conditionally on $(\xi_i)_{i=1}^N$, for all $x>0$, one has, with probability at least $1-2\exp(-x^2/2)$,
\begin{equation*}
\sup_{t \in T(\rho)} \left|\sum_{i=1}^{N} \xi_i \inr{X_i,t-t^*}\right| \leq \widehat{\sigma_N} \frac{\eta \sqrt{N} r_M^2(\rho)}{\sigma}+ x \widehat{\sigma_N} r_M(\rho)
\end{equation*} 
Thus, taking $x=\eta\sqrt{N}r_M(\rho)/\sigma$ in the previous statement, one  gets, on an event which probability is at least $1-2\exp(-\eta^2 Nr_M^2(\rho)/(2\sigma^2))$,  
\begin{equation}
\label{eq:multiplier_control2}
\sup_{t\in t^*+\rho B_1^d \cap r_M(\rho)B_2^{d}}\left|\frac{1}{N} \sum_{i=1}^{N} \xi_i \inr{X_i,t-t^*}\right| \leq 2\eta\frac{\widehat{\sigma_N}}{\sigma \sqrt{N}} r_M^2(\rho).
\end{equation}

It remains to prove, on the same event, the result for all $t\in t^*+\rho B_1^d$ such that $\norm{t-t^*}_2 > r_M(\rho)$. Define
$t^\prime:= t^* + \big(r_M(\rho)/\norm{ t-t^* }_2\big) (t-t^*)$. Since $t^\prime\in T(\rho)$, it follows from \eqref{eq:multiplier_control2} that (on the same event):
\begin{equation*}
\left|\frac{1}{N}\sum_{i=1}^{N} \xi_i \inr{X_i,\frac{r_M(\rho)}{\norm{ t-t^* }_2} (t-t^*)} \right| = \left|\frac{1}{N}\sum_{i=1}^{N} \xi_i \inr{X_i,t^\prime-t^*} \right| \leq 2\eta \frac{\widehat{\sigma_N}}{\sigma \sqrt{N}} r_M(\rho)^2
\end{equation*}
and since $r_M(\rho)\leq \norm{ t-t^* }_2$ one gets
\begin{equation*}
\left|\frac{1}{N}\sum_{i=1}^{N} \xi_i \inr{X_i,t-t^*} \right| \leq 2\eta \frac{\widehat{\sigma_N}}{\sigma \sqrt{N}} \norm{ t-t^* } r_M(\rho) 
 \leq 2\eta \frac{\widehat{\sigma_N}}{\sigma \sqrt{N}} \norm{ t-t^*}^2.
\end{equation*}


Hence, with probability (conditionally to the $\xi_i$) at least $1-2\exp \left( -\eta^2N r_M(\rho)^2)/(2\sigma^2) \right)$, the multiplier process is controlled such that
\begin{equation}
\label{eq:multiplier_control3}
\sup_{t \in t^*+\rho B_1^d} \left|P_N {\cal M}_{t-t^*}\right|\leq 4\eta \frac{\widehat{\sigma_N}}{\sigma \sqrt{N}} \max\left(r_M(\rho)^2, \norm{ t-t^*}_2^2\right).
\end{equation}

A control of the probability measure of the event $A_k$ follows by applying the previous result to $\rho=2^{k}\rho^*$ when $\eta\leq 1/(16\sqrt{2})$ together with a control of the term $\widehat{\sigma}_N$. It follows from an union bound that, conditionally to the $\xi_i$, \eqref{eq:multiplier_control3} is satisfied for all $\rho=2^{k}\rho^*, k=0, \cdots, K_0$ on an event whose probability measure is larger than
\begin{equation}
\label{eq:proba_multiplier_control1}
1-2 \sum_{k=0}^{K_0} \exp \left( -\eta^2 N r_M(2^{k} \rho^*)^2)/(2\sigma^2) \right).
\end{equation}We handle the last term below thanks to Lemma~\ref{pseudo-geometric}.

Now, we handle the random variables $\xi_1, \ldots, \xi_N$. It appears that only a control of the empirical variance term $\widehat{\sigma}_N/\sqrt{N}$ is needed to get a fully deterministic upper bound in the right-hand term of \eqref{eq:multiplier_control3}.  It follows from Bernstein inequality for subexponential variables (cf. Theorem~1.2.7 in \cite{MR3113826}) that with probability greater than $1-2\exp (-C_5 N)$,
$$\left|\frac{1}{N} \sum_{i=1}^{N} \xi_i^2 - \sigma^2\right| \leq \sigma^2,$$ 
which implies $\widehat{\sigma_N}/\sqrt{N}\leq \sqrt{2} \sigma $. Therefore, for $\eta = 1/(16\sqrt{2})$ we have 
\begin{equation}
\label{eq:proba_multiplier_control2}
\bP\left[4\eta \widehat{\sigma_N}/(\sqrt{N}\sigma)\leq 1/4\right] \geq 1- 2\exp \left( -C_5 N  \right).
\end{equation}

Binding together \eqref{eq:multiplier_control3}, \eqref{eq:proba_multiplier_control1} and \eqref{eq:proba_multiplier_control2} gives
\begin{equation*}
\bP\left[\bigcap \limits_{k=0}^{K_0} A_k\right] \geq 1- 2\exp \left( -C_5 N  \right)- 2 \sum_{k=0}^{K_0} \exp \left( -\frac{N r_M(2^{k} \rho^*)^2}{1024\sigma^2}) \right).
\end{equation*}Finally,  Lemma~\ref{pseudo-geometric} yields the following bound
\begin{equation*}
 \sum_{k=0}^{K_0} \exp \left( -\frac{Nr_M(2^{k} \rho^*)^2)}{1024\sigma^2}   \right) \leq \frac{10}{1-\exp\left(-\frac{C_M^{(1)}\sqrt{N}\rho^* }{4096\sigma}  \right)} \exp \left(- \frac{C_M^{(1)}}{C_M^{(2)}}\frac{N r_M\bigl(\rho^*\bigr)^2 }{1024\sigma^2}  \right)
 \end{equation*} and the result follows when $\rho^*\geq 4096 \log(2)  \sigma/(C_M^{(1)} \sqrt{N})$,  which implies that the denominator of the right-hand side is greater than $1/2$. 

\endproof

\subsubsection{Conclusion: construction of the event $\Omega_0$}
\label{sec:event_omega0}
We define the event
\begin{equation*}
\Omega_0 = \Omega^* \cap \bigcap \limits_{k=0}^{K_0} A_k.
\end{equation*}It follows from Proposition~\ref{prop:event_Omega^*} and Proposition~\ref{prop:event_multiplicatif} (as well as Lemma~\ref{reverse Lemma} below), that, as long as  $\rho^*\geq 4096 \log(2)  \sigma/(C_M^{(1)} \sqrt{N})$, 
\begin{equation*}
\bP[\Omega_0]\geq 1-4\exp(-C_6 N) -40 \exp \left( -\frac{C_6^\prime  N r_M(\norm{t^*}_1)^2)}{\sigma^2} \right) 
\end{equation*} where $C_6$ and $C_6^\prime$ are absolute constants.

\subsection{Deterministic part of the proof}
\label{sec:deterministic}
We first start with some few lemmas on the growth of $r_M(\cdot)$ and $r_Q(\cdot)$. We then construct a partition of $\R^d$ depending on the behavior of function $r^2(\cdot)$ (in particular, its concavity for intermediate values is an issue; we solve it thanks to a peeling argument). We then turn to the main deterministic argument showing that $\hat t$ belongs to a $\ell_1^d$-ball of radius $\rho^*$ around $t^*$. The latter holds on the event $\Omega_0$ introduced in Section~\ref{sec:event_omega0}.

\subsubsection{Two Lemmas on the growth of $r_M$ and $r_Q$}
\begin{Lemma} 
\label{rate of growth} Let $\rho>0$ and $\phi= 4(C_M^{(2)})^2/(C_M^{(1)})^2$. If $\phi\rho \leq  \rho_0  \min(1,\eta)$, then for any $\rho^\prime\geq \phi \rho$,
$$r_M^2\left(\rho^\prime\right)>2r_M^2(\rho) \textup{ and } r^2\left(\rho^\prime\right)>2r^2(\rho).$$
\end{Lemma}

\beginproof 
 Since $r_M(\cdot)$, $r_Q(\cdot)$ and $r(\cdot)$ are non-decreasing, we only have to prove the result for $\rho^\prime=\phi \rho$. Recall that $C_M^{(2)}\geq C_M^{(1)}$ (so $\phi\geq 4$). First note that if $N\geq \zeta^\prime d$ then $r_Q(\rho)=0$ and so the second claim follows from the first one since $r(\rho) = r_M(\rho)$ in this case. And when $N\leq \zeta d$, one has $r_Q^2(\phi\rho)=\phi^2r_Q(\rho)\geq 16 r_Q(\rho)> 2 r_Q(\rho)$ (because in that  case $r_Q(\rho)^2>0$). Therefore the second claim is a straightforward consequence of the first one. So it only remains to study the behavior of $r_M^2(\cdot)$. For $\rho < \phi \rho \leq \rho_0$, $r_M$ is given by one of the two last expressions of~\eqref{eq:rM}. 

First assume that $(\phi \rho)^2 N \leq  \sigma^2 \log d $ then 
\begin{equation*}
r_M^2(\phi \rho) \geq C_M^{(1)}\phi \rho \sigma \sqrt{\log(ed/N)}>2C_M^{(2)}\rho \sigma \sqrt{\log(ed/N)}\geq 2r_M^2(\rho).
\end{equation*}  

Now, assume that $\sigma^2 \log d \leq \rho^2 N   \leq (\phi \rho) ^2 N$ then 
$$r_M^2(\phi \rho) = C_M  \phi \rho \sigma \sqrt{\frac{1}{N}\log\left(\frac{e\sigma d}{\phi \rho\sqrt{N}}\right)}$$ for some $C_M\in [C_M^{(1)},C_M^{(2)}]$. One has that for all $x\geq \phi e$, $\log(x/\phi)>\log(x)/\phi$, and, since $\rho_0=\sigma d/(\eta \sqrt{N})$,the assumption $\phi\rho \leq  \rho_0 \, \min(1,\eta)$ guarantees that $e\sigma d/\rho \sqrt{N} \geq  \phi e$. Therefore, we have:
\begin{equation*}
r_M^2(\phi \rho) > C_M^{(1)} \phi\rho \sigma \sqrt{\frac{1}{\phi N}\log\Big(\frac{e\sigma d}{\rho\sqrt{N}}\Big)} \geq 2C_M^{(2)}\rho \sigma \sqrt{\frac{1}{N}\log\Big(\frac{e\sigma d}{\rho\sqrt{N}}\Big)} \geq 2 r_M^2(\rho).
\end{equation*}

Finally, when  $ \rho^2 N \leq \sigma^2 \log d    \leq (\phi \rho) ^2 N$, since $r_M$ is increasing, it is clear considering the two previous cases that one has again $r_M^2(\phi \rho) > (\sqrt{\phi}C_M^{(1)} / C_M^{(2)}) r_M^2(\rho)$, so $r_M^2(\phi \rho) >2 r_M^2(\rho)$.

\endproof

\begin{Lemma}
\label{reverse Lemma} Let $\nu>0$. 
If $\nu \geq 1$ then $r_M(\nu \rho)\leq \sqrt{\nu}r_M(\rho)$ and  $r(\nu \rho)\leq \nu r(\rho)$. If $\nu \leq 1$ then $r_M(\nu \rho)\geq \sqrt{\nu}r_M(\rho)$ and $r(\nu \rho)\geq \nu r(\rho)$.
\end{Lemma}
 
 \beginproof
It is clear that $r_Q(\nu \rho)=\nu r_Q(\rho)$, because $\ell^*(\nu \rho B_1^d \cap \nu r_Q(\rho) B_2^d)=\nu \ell^*(\rho B_1^d \cap r_Q(\rho) B_2^d)$. As for $r_M(\nu \rho)$, for $\nu \geq 1$, one has 
\newline $\sigma \ell^*(\nu \rho B_1^d \cap \sqrt{\nu} r_M(\rho) B_2^d)\leq \sigma \nu \ell^*(\rho B_1^d \cap r_M(\rho) B_2^d)= \nu \eta \sqrt{N} r_M(\rho)^2= \eta \sqrt{N} (\sqrt{\nu}r_M(\rho))^2$. So $r_M(\nu \rho)\leq \sqrt{\nu}r_M(\rho)$.

As for the case $\nu\leq 1$, then $1/\nu\geq 1$, and it suffices to write that $r_M(\rho)=r_M\left( (1/\nu)\nu \rho \right)\leq (1/\nu) r_M(\nu \rho)$ to get the result (and still $r_Q(\rho/\nu)=r_Q(\rho)/\nu$).

 \endproof

\subsubsection{Partition of $\R^d$ into three zones}
\label{sec:zones}
Recall that $\rho^*=\max\left(10,8(C_M^{(2)})^2/((C_M^{(1)})^2\eta)+1\right) \norm{t^*}_1$, that $\rho_0$ is the smallest $\rho$ such that  $r_M(\rho)=r_M(\rho_0)$ for all $\rho\geq \rho_0$ and that  $K_0=\min \{k \in \mathbb{N}: 2^k \rho^* \geq 2\rho_0\}$.

Hereunder, $K_0$ is the number of shelves used to partition the intermediate ``peeling zone'' defined below. We use $\rho^*$ and $K_0$ to construct a  partition of $\mathbb{R}^d$ into three main zones: 
\begin{itemize}
\item the ``central zone'' $t^*+\rho^*B_1^d$,
\item the intermediate ``peeling zone'' $\{t\in\R^d: \rho^*<\norm{t-t^*}_1\leq 2^{K_0}\rho^*\}$ (recall that $2^{K_0}\rho^* \simeq 2\rho_0$). This zone is considered only when $K_0\geq 1$. We use a ``peeling", i.e. a partition of this zone into $K_0$ sub-areas called the ``shelves":  $2^{k-1} \rho^*<\norm{t-t^*}_1\leq 2^{k} \rho^*$, for $k=1, \cdots, K_0$.

\item the ``exterior zone'' $\{t\in\R^d: \norm{t-t^*}_1> 2^{K_0}\rho^*\}$, on which  $r_M(\cdot)$ is constant.
\end{itemize}

Our main objective is now to show that, on the event $\Omega_0$, $\hat{t}$ belongs to the central zone.

\subsubsection{Locating $\hat{t}$ in the central zone on the event $\Omega_0$}
To show that, on the event $\Omega_0$, $\hat{t}$ belongs to the central zone, it is enough to show that any $t$ outside this area satisfies $P_N\cL_t^\Psi>0$ where we recall that in our case,
\begin{equation*}
P_N\cL_t^\Psi:=\left(\frac{1}{N}\sum_{i=1}^N (Y_i-\inr{X_i,t})^2 + c_0 r^2(\norm{t}_1)\right)- \left(\frac{1}{N}\sum_{i=1}^N (Y_i-\inr{X_i,t^*})^2 + c_0 r^2(\norm{t^*}_1)\right).
\end{equation*} This will indeed prove that $\norm{\hat t-t^*}_1\leq \rho^*$ since $P_N\cL_{\hat{t}}^\Psi\leq 0=P_N\cL_{t^*}^\Psi$.

Let $t\in\R^d$ be outside the central zone -- i.e. $\norm{t-t^*}_1>\rho^*$. First note that $t$ satisfies $\norm{t}_1 \geq\norm{t-t^*}_1-\norm{t^*}_1 
\geq \rho^*-\norm{t^*}_1 \geq 9\norm{t^*}_1$. Therefore, we have ${\cal R}_{t^\prime,t^*} = \Psi(\norm{t}_1)-\Psi(\norm{t^*}_1) = c_0 r^2(\norm{t}_1) - c_0 r^2(\norm{t^*}_1)>0$ and   $(8/9)\norm{t}_1 \leq \norm{t-t^*}_1  \leq (10/9)\norm{t}_1$.

\textbf{From now on, we place ourselves on the event $\Omega_0$ introduced in Section~\ref{sec:event_omega0}.}

First suppose that $\rho^*<2\rho_0$,  then both the ``intermediate peeling zone'' and the ``exterior zone'' must be considered.  For $t$ in any of these two areas, one has ${\cal R}_{t^\prime,t^*}\geq 0$.

Let us begin by the ``peeling zone". Consider $t$ and $k\geq 1$ such that
\begin{equation*}
2^{k-1}\rho^*<\norm{t-t^*}_1\leq 2^{k}\rho^* \textup{ and } \norm{t-t^*}_1 \leq 2^{K_0}\rho^*.
\end{equation*} We recall that $2^{K_0-1}\rho^*<2\rho_0\leq 2 ^{K_0}\rho^*$ and $\rho_0$ is the smallest radius such that $r_M(\rho) = r_M(\rho_0)$ for all $\rho\geq \rho_0$.

One possibility is that  $\norm{t-t^*}_2^2\geq r^2(2^{k}\rho^*)$, then on $\Omega^*$ one has $P_N {\cal Q}_{t-t^*}\geq \norm{t-t^*}_2^2/2$ and, on $A_k$, one has $|P_N {\cal M}_{t-t^*}|\leq  \norm{t-t^*}_2^2/4$. Hence, on $\Omega_0$,
$P_N\cL_t^\Psi \geq \norm{t-t^*}_2^2/4+{\cal R}_{t^\prime,t^*} 
>0.$

Let us tackle now the other case: $\norm{t-t^*}_2^2<r^2(2^{k}\rho^*)$. We will show that in this situation, if $c_0$ is large enough, ${\cal R}_{t,t^*}/4>|P_N {\cal M}_{t-t^*}|$.

As $\norm{t}_1 > 4\norm{t^*}_1(C_M^{(2)})^2/(C_M^{(1)})^2$ ($t$ is not in the ``central zone'') and $4\norm{t^*}_1(C_M^{(2)})^2/(C_M^{(1)})^2\leq \eta \rho_0$ (since $2\rho_0\geq \rho^*> 8\norm{t^*}_1(C_M^{(2)})^2/((C_M^{(1)})^2 \eta)$), by Lemma~\ref{rate of growth} one has that 
$c_0 r^2(\norm{t}_1)\geq 2c_0 r^2(\norm{t^*}_1).$ Thus, ${\cal R}_{t,t^*} \geq c_0 r^2(\norm{t}_1)/2$. So, thanks to Lemma~\ref{reverse Lemma}, and since 
$$\frac{\norm{t}_1}{2^{k}\rho^*}=\frac{\norm{t-t^*}_1}{2^{k}\rho^*}\frac{\norm{t}_1}{\norm{t-t^*}_1}>\frac{1}{2}\cdot\frac{9}{10}> \frac{2}{5},$$ 
 one has 
 $${\cal R}_{t,t^*} \geq \frac{1}{2}c_0 r^2(\norm{t}_1)\geq \frac{1}{2}\cdot\left(\frac{\norm{t}_1}{2^{k}\rho^*}\right)^2 c_0 r^2(2^{k}\rho^*) >(2/25)c_0r^2(2^{k}\rho^*).$$
In addition, we have $\norm{t-t^*}_1\leq 2^{k}\rho^*$ and $\norm{t-t^*}_2^2<r^2(2^{k}\rho^*)$, hence, on the event $A_{k}$, $|P_N {\cal M}_{t-t^*}|\leq r^2(2^{k}\rho^*)/4$. As a consequence, for $c_0>13$ one has 
\begin{equation}
\label{eq:marge1}
 \frac{{\cal R}_{t,t^*}}{4}  -|P_N {\cal M}_{t-t^*}| >\frac{2}{25}\cdot c_0r^2(2^{k}\rho^*)\cdot\frac{1}{4}-\frac{1}{4}r^2(2^{k}\rho^*)>0.
 \end{equation}
A fortiori, ${\cal R}_{t^\prime,t^*}-|P_N {\cal M}_{t-t^*}| >0$ and as $P_N {\cal Q}_{t-t^*}\geq0$, this implies in particular that  $P_N\cL_t^\Psi>0$.

To sum up, we proved that for all $t$ in the peeling zone, we have
\begin{equation}
\label{eq:peeling}
P_N\cL_t^\Psi\geq P_N {\cal Q}_{t-t^*}-|P_N {\cal M}_{t-t^*}|+{\cal R}_{t^\prime,t^*}/4>0.
\end{equation} 

Let us now study the exterior zone. We will mainly use homogeneity arguments. 
Let $t\in\R^d$ be outside the ball $t^*+2^{K_0}\rho^*B_1^d$. Define $t^\prime\in t^*+2^{K_0}\rho^* S_1^{d-1}$ such that  $t-t^*=\alpha_t (t^\prime-t^*)$ for some $\alpha_t\geq1$. In particular, $\norm{t^\prime-t^*}_2 \geq 2 r_0> r_M(2^{K_0}\rho^*)$ so by Proposition~\ref{prop:event_multiplicatif}, on the event $A_{K_0}$, one has $|P_N {\cal M}_{t^\prime-t^*}|\leq \norm{t^\prime-t^*}_2^2/4$. We consider now two cases. 

If $\norm{t^\prime-t^*}_2\geq  r_Q(\norm{t^\prime-t^*}_1)$, then by Lemma~\ref{prop:event_Omega^*}, on the event $\Omega^*$, one has $P_N {\cal Q}_{t^\prime-t^*}\geq \norm{t^\prime-t^*}_2^2/2$. So on the event $\Omega_0$, $P_N {\cal Q}_{t^\prime-t^*}-|P_N {\cal M}_{t^\prime-t^*}|>0$ and 
\begin{equation*}
P_N {\cal Q}_{t-t^*}-|P_N {\cal M}_{t-t^*}|=\alpha_t^2 P_N {\cal Q}_{t^\prime-t^*}-\alpha_t|P_N {\cal M}_{t^\prime-t^*}|\geq \alpha_t\left(P_N {\cal Q}_{t^\prime-t^*}-|P_N {\cal M}_{t^\prime-t^*}|\right)>0
\end{equation*} therefore $P_N\cL_t^\Psi >0$.

On the contrary, if $\norm{t^\prime-t^*}_2<r_Q(\norm{t^\prime-t^*}_1)$, as $\norm{t^\prime-t^*}_2\geq 2r_0$, then $r_Q(\norm{t^\prime}_1)\geq 9r_Q(\norm{t^\prime-t^*}_1)/10>9r_0/5= 9r_M(\norm{t^\prime}_1)/5>r_M(\norm{t^\prime}_1) $ so $r(\norm{t^\prime}_1)=r_Q(\norm{t^\prime}_1)$. One has $\norm{t}_1=\norm{t^*+\alpha_t (t^\prime-t^*)}_1 
 \geq (4\alpha_t/5)\norm{t^\prime}_1$ since $\alpha_t \geq 1$ and 
$\norm{t^\prime-t^*}_1\geq \rho^*$. So  $r^2(\norm{t}_1) \geq r_Q^2\left(\alpha_t4\norm{t^\prime}_1/5\right)=4^2\alpha_t^2r^2(\norm{t^\prime}_1)/5^2.$ We have seen before that $r_Q(\norm{t^\prime}_1)\geq 9r_0/5 \geq 9r_M(\norm{t^*}_1)/5$, and $r_Q(\norm{t^\prime}_1)\geq 9 r_Q(\norm{t^*}_1)$, so $r^2(\norm{t^\prime}_1)\geq 3 r^2(\norm{t^*}_1)$. As a consequence, 
\begin{align*} 
{\cal R}_{t,t^*} & = c_0r^2(\norm{t}_1)-c_0r^2(\norm{t^*}_1)  \geq c_0 \left( \alpha_t^2\frac{4^2}{5^2}r^2(\norm{t^\prime}_1)-(1/3)r^2(\norm{t^\prime}_1)\right) \\
&  \geq \alpha_t \frac{c_0}{4} r^2(\norm{t^\prime}_1)  \geq \alpha_t {\cal R}_{t^\prime,t^*}/4.
\end{align*}
 Moreover, since $t-t^*=\alpha_t (t^\prime-t^*)$, one has $P_N {\cal Q}_{t-t^*}=\alpha_t^2 P_N {\cal Q}_{t^\prime-t^*}$ and $P_N {\cal M}_{t-t^*}=\alpha_t P_N {\cal M}_{t^\prime-t^*}$. So, in the case $\norm{t^\prime-t^*}_2<r_Q(\norm{t^\prime-t^*}_1)$, by \eqref{eq:peeling} applied to $t^\prime$,  
$$P_N\cL_t^\Psi \geq \alpha_t (P_N {\cal Q}_{t^\prime-t^*}-|P_N {\cal M}_{t^\prime-t^*}|+ {\cal R}_{t^\prime,t^*}/4)>0.$$

Let us now consider the case $\rho^*>2\rho_0$. In this situation,  there is no need for the intermediate ``peeling" zone. Let $t\in\R^d$ be outside the ball  $t^*+\rho^*B_1^d$ and  set $t^\prime\in t^*+\rho^*S_1^{d-1}$ such that  $t-t^*=\alpha_t (t^\prime-t^*) \textup{ with } \alpha_t\geq1$.  Then one can apply arguments similar to the peeling case (with $\norm{t^\prime-t^*}_1=\rho^*$, but this time $k=K_0=0$), on the event $A_0\cap \Omega^*$, to $t^\prime$. If $\norm{t^\prime-t^*}_2\geq  r_Q(\norm{t^\prime-t^*}_1)$, then by Lemma~\ref{prop:event_Omega^*}, on the event $\Omega_0$ one has $P_N {\cal Q}_{t^\prime-t^*}\geq |P_N {\cal M}_{t^\prime-t^*}|/2$. Conversely, if $\norm{t^\prime-t^*}_2< r_Q(\norm{t^\prime-t^*}_1)$, then $r^2(\norm{t^\prime}_1)\geq 3 r^2(\norm{t^*}_1)$ so for $c_0$ big enough, in the same spirit as~\eqref{eq:marge1}, one gets that ${\cal R}_{t^\prime,t^*}/4- |P_N {\cal M}_{t^\prime-t^*}|>0$. So in both cases, $P_N {\cal Q}_{t-t^*}-|P_N {\cal M}_{t-t^*}|+{\cal R}_{t^\prime,t^*}/4>0$. The same argument as previously for the exterior zone, shows that $P_N\cL_t^\Psi\geq \alpha_t (P_N {\cal Q}_{t^\prime-t^*}-|P_N {\cal M}_{t^\prime-t^*}|+{\cal R}_{t^\prime,t^*}/4)$, so $P_N\cL_t^\Psi>0$.
As a conclusion, in the case $\rho^*>2\rho_0$, we have for all $t\in\R^d$ satisfying $\norm{t-t^*}_1\geq \rho^*>2\rho$ that $P_N\cL_t^\Psi>0$.

To sum up, on the event $\Omega_0$, any $t$ outside the central zone satisfies $P_N\cL_t^\Psi>0$. Therefore, given that $\hat{t}$ satisfies $P_N\cL_{\hat t}^\Psi\leq 0$, we conclude that $\hat t$ belongs to the central zone.

\subsubsection{Conclusion of the proof of Theorem~\ref{theo:main}}
\label{end proof}
On the event $\Omega_0$, $\hat t\in (t^*+\rho^* B_1^d)$. Hence, either $\norm{\hat{t}-t^*}_2^2\leq r^2(\rho^*)$ and the proof is over or $\norm{\hat{t}-t^*}_2^2>r^2(\rho^*)$. In the latter case, one has
\begin{equation*}
P_N {\cal Q}_{\hat t-t^*}>\frac{1}{2} \norm{\hat t-t^*}_2^2 \textup{ and }|P_N {\cal M}_{\hat t-t^*}|<\frac{1}{4} \norm{\hat t-t^*}_2^2
\end{equation*} and so
$$ 0 \geq P_N\cL_{\hat{t}}^\Psi \geq \frac{1}{4} \norm{\hat{t}-t^*}_2^2 +c_0r^2(\norm{\hat{t}}_1)-c_0r^2(\norm{t^*}_1) $$
which implies  $ \norm{\hat{t}-t^*}_2^2 \leq 4c_0r^2(\norm{t^*}_1) .$  

Thus, taking $\theta_0= \max \left(100,(8(C_M^{(2)})^2/(C_M^{(1)})^2+1)^2, 4c_0\right),$ with probability at least 
$$\bP[\Omega_0] \geq 1-4\exp \left( -C_6 N \right) -40 \exp \left( - C_6^\prime N r_M(\norm{t^*}_1)^2/\sigma^2\right)$$
one gets that in both cases, $R(\hat t) - R(t^*)=\norm{\hat{t}-t^*}_2^2 \leq \theta_0 r^2(\norm{t^*}_1)$. Moreover, for $\norm{t^*}_1\geq \Delta_0\sigma \sqrt{\log(ed)/N}$ for $\Delta_0$ an absolute constant large enough, we have $\bP[\Omega_0]\geq 3/4$. Given that $r^2(\norm{t^*}_1)$ is the minimax rate of convergence over $\norm{t^*}_1 B_1^d$ (cf. Theorem~\ref{thm:minimax_B1}), we conclude that $\Psi(\rho) = c_0 r^2(\rho)$
is indeed a minimax regularization function.

\section{Technical material and proof of Proposition~\ref{prop:adapt_small_ball}}
\label{sec:technical_material}

\subsection{Localization with balls and spheres}

The next lemma shows that when the intersection is not trivial, localizing by intersecting an $\ell_1$ ball with an $\ell_2$ sphere or the corresponding full $\ell_2$ ball is equivalent. 

\begin{Lemma}
\label{prop: equivalent localizations}
If $\rho\geq r$, then $\ell^*(\rho B_1^d \cap r S_2^{d-1})=\ell^*(\rho B_1^d \cap r B_2^d).$ If $\rho<r$, then $\ell^*(\rho B_1^d \cap r B_2^d)=\ell^*(\rho B_1^d)$ and $\ell^*(\rho B_1^d \cap r S_2^{d-1})=0$.
\end{Lemma}

\beginproof
Since for any set $T\subset \R^d$, $\ell^*(T)=\ell^*({\rm conv}(T))$ (with ${\rm conv}(T)$ denoting the convex hull of $T$), the result is a direct consequence of the fact that for $\rho<r$, ${\rm conv}(\rho B_1^d \cap r B_2^d)=\rho B_1^d$ and ${\rm conv}(\rho B_1^d \cap r S_2^{d-1})=\emptyset$ (which are two obvious statements), and that if $\rho\geq r$, then ${\rm conv}(\rho B_1^d \cap r S_2^{d-1})={\rm conv}(\rho B_1^d \cap r B_2^d)$, which we prove now. 

One inclusion is immediate, it remains to prove that ${\rm conv}(\rho B_1^d \cap r B_2^d)\subset {\rm conv}(\rho B_1^d \cap r S_2^{d-1})$. First, 
$${\rm conv}(\rho B_1^d \cap r B_2^d)=\rho B_1^d \cap r B_2^d={\rm conv}\{t\in\R^d: \ \max(\norm{t}_1/\rho, \norm{t}_2/r)=1\}$$
 so it  only remains to show that $\{t\in\R^d: \ \max(\norm{t}_1/\rho, \norm{t}_2/r)=1\}\subset {\rm conv}(\rho B_1^d \cap r S_2^{d-1})$.
 
 First, remark that $\{t\in\R^d: \norm{t}_1/\rho\leq 1, \norm{t}_2/r=1\}$ is included in $\rho B_1^d \cap r S_2^{d-1}$.  Let us now consider the set $\{t\in\R^d: \norm{t}_1/\rho= 1, \norm{t}_2/r<1\}$ and consider an element $t$ in it. 
 We denote by $e_1,\ldots, e_d$ the canonical basis of $\R^d$ and we recall that each face of $\rho B_1^d$ is the convex hull of its vertices. So, since $t \in \rho S_1^{d-1}$, there exist $b_1\in \{\rho e_1,-\rho e_1\}, b_2 \in \{\rho e_2,-\rho e_2\},\ldots, b_d\in \{\rho e_d,-\rho e_d\}$ such that $t$ is in the convex hull of $\{ b_1\ldots,  b_d\}$: there exist non-negative coefficients $\mu_1,\ldots, \mu_d$ such that $\sum_{j=1}^{d} \mu_j  (b_j-t)=0$. Consider for each $j \in \{1,\ldots,d\}$ the mapping 
 $$f_j: x\in [0,1] \mapsto \norm{t+x(b_j-t)}_2/r-\norm{t+x(b_j-t)}_1/\rho$$
 This mapping is continuous, $f_j(0)< 0$ ($t$ is in $\rho S_1^{d-1}$ and $r B_2^d$ but not in $r S_2^{d-1}$) , and $f_j(1)\geq 0$ (because $\norm{b_j}_2=\norm{b_j}_1=\rho$ and $\rho\geq r$). So there exists $x_j$ in (0,1] such that $f(x_j)=0$, which means that $t+x_j(b_j-t) \in \rho S_1^{d-1} \cap r S_2^{d-1}$. One has  that
 $$\sum \limits_{j=1}^{d} \frac{\mu_j}{x_j} \left(t+x_j  (b_i-t)-t\right)=\sum \limits_{j=1}^{d} \frac{\mu_j}{x_j} x_j  (b_j-t)=\sum \limits_{j=1}^{d} \mu_j  (b_j-t)=0.$$
  As a consequence, $t$ is in the convex hull of the vectors $t+x_j  (b_j-t) \in \rho S_1^{d-1}\cap r S_2^{d-1}, j\in\{1, \ldots, d\}$ which achieves the proof.
\endproof

\subsection{Control of the probability estimate}

\begin{Lemma}
\label{pseudo-geometric}Set $\eta=1/(16\sqrt{2})$. For all $\nu>0$, we have
\begin{equation*}
\sum_{k=0}^{K_0} \exp \left(-\nu r_M \bigl(2^{k}\rho^*\bigr)^2\right) \leq \frac{10}{1-\exp\left(-\nu C_M^{(1)} \frac{\sigma}{4\sqrt{N}} \rho^* \right)} \exp \left(- \frac{C_M^{(1)}}{C_M^{(2)}}\nu  r_M\bigl(\rho^*\bigr)^2 \right).
\end{equation*}
\end{Lemma}

\beginproof
First, the terms of the sum are non-increasing (remember that $r_M$ is a non-decreasing function). So skipping the last terms will not change the order of magnitude: $$\sum_{k=0}^{K_0} \exp \left(-\nu r_M \bigl(2^{k}\rho^*\bigr)^2\right) \leq \max\left(10\exp \left(-\nu r_M \bigl(\rho^*\bigr)^2\right),10 \sum_{k=0}^{K_0-9} \exp \left(-\nu r_M \bigl(2^{k}\rho^*\bigr)^2\right)\right).$$

One has $10\exp \left(-\nu r_M \bigl(\rho^*\bigr)^2\right) > 10 \sum_{k=0}^{K_0-9} \exp \left(-\nu r_M \bigl(2^{k}\rho^*\bigr)^2\right)$ when $K_0\leq 8$. Let us now assume that $K_0\geq 9$. We study the sum $\sum_{k=0}^{K_0-9}\exp \left(-\nu r_M \bigl(2^{k}\rho^*\bigr)^2\right)$.

 In order to get rid of the ``range" $[C_M^{(1)},C_M^{(2)}]$ in the definition of $r_M(\cdot)$, notice that $\sum_{k=0}^{K_0-9} \exp \left(-\nu r_M \bigl(2^{k}\rho^*\bigr)^2\right)  \leq \sum_{k=0}^{K_0-9} a_k $ with:
\begin{equation*}
 a_k:=\left\{ 
 \begin{array}{cc}
 \exp \left(- \nu  C_M^{(1)} 2^{k}\rho^* \sigma \sqrt{\frac{\log(ed)}{N}}\right) & \mbox{ if } \bigl(2^{k}\rho^*\bigr)^2 N\leq \sigma^2 \log(d)
\\
\exp \left(- \nu C_M^{(1)} 2^{k}\rho^* \sigma\sqrt{\frac{1}{N}\log\Big(\frac{e\sigma d}{2^{k}\rho^*\sqrt{N}}\Big)} \right) & \mbox{otherwise.}
 \end{array}
\right.
\end{equation*}

We emphasize that the sum goes only up to $2^{K_0} \rho^*$, which excludes the constant third form of $r_M(\cdot)$.

Applying a second time the ``range" $[C_M^{(1)},C_M^{(2)}]$ in the bounds on $r_M$ allows to bound $a_0$ in terms of $r_M\bigl(\rho^*\bigr)^2$: $a_0 \leq \exp \left(- (C_M^{(1)}/C_M^{(2)})\nu  r_M\bigl(\rho^*\bigr)^2 \right).$
Therefore, in the following we will bound $\sum_{k=0}^{K_0-9} a_k $ with respect to $a_0$, and then get back to $r_M \bigl(\rho^*\bigr)^2$.

We now prove that there exists $\alpha$ independent on $k$ (but dependent on the other parameters) such that for any $k\leq K_0-9$,  $a_{k+1}/a_k\leq \alpha <1$.  

If $(2^{k+1}\rho^*)^2\leq  \sigma^2 \log(d)/N$ then
\begin{align*}
\frac{a_{k+1}}{a_k} &=\exp \left(-\nu C_M^{(1)} 2^{k+1}\rho^* \sigma \sqrt{\frac{\log(e d)}{N}}+ \nu C_M^{(1)} 2^{k}\rho^* \sigma \sqrt{\frac{\log(e d)}{N}} \right)  \\
&=\exp \left(-\nu C_M^{(1)} 2^{k}\rho^* \sigma \sqrt{\frac{\log(e d)}{N}} \right) \leq \exp(-\nu C_M^{(1)} \rho^* \beta_1)
 \end{align*}
 with $\beta_1=\sigma \sqrt{\log(e d)/N}\geq \sigma/\sqrt{N}$.
 
 As for the case $(2^k\rho^*)^2 \leq \sigma^2 \log(d)/N < (2^{k+1}\rho^*)^2 $ (which can only occur when $d>1$), then one has $2^{k+1}\rho^* \leq 2 \sigma \sqrt{\log d}/\sqrt{N}$ and so
\begin{align*}
\frac{a_{k+1}}{a_k} &= \exp \left(- \nu C_M^{(1)} 2^{k+1}\rho^* \sigma\sqrt{\frac{1}{N}\log\Big(\frac{e\sigma d}{2^{k+1}\rho^*\sqrt{N}}\Big)}+\nu  C_M^{(1)} 2^{k}\rho^* \sigma \sqrt{\frac{\log(ed)}{N}} \right) \\
&\leq \exp \left(- \nu C_M^{(1)} 2^k \rho^* \frac{\sigma}{\sqrt{N}} \left( 2\sqrt{\log \left( \frac{ed}{2\sqrt{\log d}} \right)}- \sqrt{\log ed}\right) \right) \\
& \leq \exp \left(- \nu C_M^{(1)} 2^k \rho^* \frac{\sigma}{2\sqrt{N}} \right)  \leq \exp \left(- \nu C_M^{(1)}\rho^* \frac{\sigma}{2\sqrt{N}}\right)
  \end{align*}
The second inequality relies on a straightforward analysis fact:
$$\forall \ d \geq 2, \quad \left( 2\sqrt{\log \left( \frac{ed}{2\sqrt{\log d}} \right)}- \sqrt{\log ed}\right)>\frac{1}{2} $$
 
 Let us tackle now the case $(2^k\rho^*)^2>\sigma^2 \log(d)/N$. We have $a_{k+1}/a_k=\exp(b_k \nu C_M^{(1)} 2^k\rho^* \sigma/\sqrt{N})$ where
$$b_k:=-2\sqrt{\log\left(\frac{e\sigma d}{2^{k+1}\rho^*\sqrt{N}}\right)}+\sqrt{\log\left(\frac{e\sigma d}{2^{k+1}\rho^*\sqrt{N}}\right)+ \log(2)}.$$ Since $\sqrt{x+y}\leq \sqrt{x}+\sqrt{y}$ for all $x,y\geq0$, one has (still for $k\leq K_0-9$): 
$$b_k\leq \sqrt{\log(2)}-\sqrt{\log\left(\frac{e\sigma d}{2^{k+1}\rho^*\sqrt{N}}\right)} \leq \sqrt{\log(2)}-\sqrt{\log\left(\frac{e\sigma d}{2^{K_0-8}\rho^*\sqrt{N}}\right)}\leq (1-\sqrt{2}) \sqrt{\log(2)}$$
because $2^{K_0-8}\rho^*\leq 2^{-6} \rho_0 \leq 2^{-6} \sigma d/(\eta \sqrt{N})\leq \sigma d/(2\sqrt{N})$ (we recall that $\rho_0=\sigma d/(\eta \sqrt{N})$ and $\eta=1/(16\sqrt{2})$). It follows that 
$$\frac{a_{k+1}}{a_k} =\exp \left(\nu C_M^{(1)} 2^k \rho^* \frac{\sigma}{\sqrt{N}} b_k \right) \leq \exp \left(-\nu C_M^{(1)} \rho^* \frac{\sigma}{\sqrt{N}} (-b_k) \right)\leq \exp(-\nu C_M^{(1)} \rho^* \beta_2)$$ with $\beta_2=(\sqrt{2}-1)\sqrt{\log(2)}\sigma/\sqrt{N}\geq \sigma/(4\sqrt{N})$. Then, we conclude that for all  $k \leq K_0-9$, $a_{k+1}/a_k\leq \exp\left(-\nu C_M^{(1)} \frac{\sigma}{4\sqrt{N}} \rho^* \right)$ and so
$$\sum_{k=0}^{K_0-9} a_k \leq a_0 \sum_{k=0}^{K_0-9}\exp\left(-k\nu C_M^{(1)} \frac{\sigma}{4\sqrt{N}} \rho^* \right)=a_0 \frac{1-\exp\left(-(K_0-8)\nu C_M^{(1)} \frac{\sigma}{4\sqrt{N}} \rho^* \right)}{1-\exp\left(-\nu C_M^{(1)} \frac{\sigma}{4\sqrt{N}} \rho^* \right)}.$$
Finally, the result follows, since $a_0 \leq \exp \left(- (C_M^{(1)}/C_M^{(2)})\nu  r_M\bigl(\rho^*\bigr)^2 \right)$ and 
$$ \frac{1-\exp\left(-(K_0-8)\nu C_M^{(1)} \frac{\sigma}{4\sqrt{N}} \rho^* \right)}{1-\exp\left(-\nu C_M^{(1)} \frac{\sigma}{4\sqrt{N}} \rho^* \right)} \leq \frac{1}{1-\exp\left(-\nu C_M^{(1)} \frac{\sigma}{4\sqrt{N}} \rho^* \right)}.$$
\endproof

\subsection{Proof of Proposition~\ref{prop:adapt_small_ball}} 
\label{sub:proof_of_prop}
The proof of Proposition~\ref{prop:adapt_small_ball} relies on key ideas developed in minimax theory. We refer the reader to \cite{MR2724359} for a state of the art in minimax theory. 

Let $\hat t$ satisfy the properties of Proposition~\ref{prop:adapt_small_ball}, i.e. a procedure adaptive on the finite set $\{\pm \rho e_1, \cdots, \pm \rho e_d\}$ and let $t^*\in (\rho/2)B_1^d$.  Denote $\Lambda = \{t^*, \pm \rho e_1, \ldots, \pm \rho e_d\} = \{t_0^*, t_1^*, \cdots, t_{2d}^*\}$ so that $t_0^*=t^*$. It is straightforward to check that $\Lambda$ is a $\rho/2$-separated set in $\rho B_1^d$ w.r.t. $\ell_2^d$. 



Now, let us define the test statistics 
 \begin{equation*}
 \hat \phi \in\argmin_{j\in\{0, \cdots, 2d\}} \norm{t^*_j - \hat t}_2.
 \end{equation*}One has for all $j\in\{0, 1, \ldots, 2d\}$ that, if $\hat \phi\neq j$ then $\norm{\hat t -t^*_j}_2\geq \rho/4$. Indeed, if $\hat \phi\neq j$ then there exists $k\in\{0, 1, \dots, 2d\}\backslash\{j\}$ such that $\norm{t^*_k-\hat t}_2\leq \norm{t^*_j-\hat t}_2$. If $\norm{t^*_k-\hat t}_2\geq \rho/4$ then the result holds otherwise $\norm{t^*_k-\hat t}_2< \rho/4$ and so $\norm{t^*_j-\hat t}_2\geq \norm{t^*_j- t^*_k}_2-\norm{t^*_k-\hat t}_2\geq \rho/4$. Therefore, we have, for all $\tau>0$
 \begin{align}\label{eq:minimax_small_ball_main}
 \nonumber &\bP_{t_0^*}\left[\norm{\hat t - t_0^*}_2\geq \rho/4\right]\geq \bP_{t_0^*}\left[\hat\phi \neq0\right] = \sum_{j=1}^{2d} \bP_{t_0^*}\left[\hat\phi = j\right]\geq \sum_{j=1}^{2d} \tau \bP_{t_j^*}\left[ \hat\phi = j \mbox{ and } \frac{d\bP_{t_0^*}}{d\bP_{t_j^*}}\geq \tau\right]\\
 &\geq \tau \sum_{j=1}^{2d}\bP_{t_j^*}\left[ \hat\phi = j\right] - \bP_{t_j^*}\left[ \frac{d\bP_{t_0^*}}{d\bP_{t_j^*}}< \tau\right]\geq \tau \sum_{j=1}^{2d}\bP_{t_j^*}\left[ \norm{\hat t-t_j^*}_2<\rho/4\right] - \bP_{t_j^*}\left[ \frac{d\bP_{t_0^*}}{d\bP_{t_j^*}}< \tau\right].
 \end{align}

 It follows from the adaptation property of $\hat t$ over $\{\pm \rho e_1, \cdots, \pm \rho e_d\}$ that for every $j\in\{1, \cdots, 2d\}$, 
 \begin{equation}\label{eq:main_adapt}
 \bP_{t_j^*}\left[ \norm{\hat t-t_j^*}_2<\rho/4\right] \geq \bP_{t_j^*}\left[ \norm{\hat t-t_j^*}_2^2< \chi_1 r^2(\rho)\right]\geq 3/4 
 \end{equation}when $\chi_1 r^2(\rho)\leq \rho^2/16$.
 Let $j\in\{1, \cdots, 2d\}$. Following the same argument as in Proposition~2.3 in \cite{MR2724359}  (based on second Pinsker inequality), we obtain
 \begin{align*}
 &\bP_{t_j^*}\left[ \frac{d\bP_{t_0^*}}{d\bP_{t_j^*}}\geq  \tau\right] = \bP_{t_j^*}\left[ \frac{d\bP_{t_j^*}}{d\bP_{t_0^*}}\leq \frac{1}{\tau}\right] = 1 - \bP_{t_j^*}\left[ \log\frac{d\bP_{t_j^*}}{d\bP_{t_0^*}}> \log(1/\tau)\right]\\
 &\geq 1- \frac{1}{\log(1/\tau)}\int \left(\log\frac{d\bP_{t_j^*}}{d\bP_{t_0^*}}\right)_+ d \bP_{t^*_j}\geq 1-\frac{1}{\log(1/\tau)}\left[K(\bP_{t_j^*}, \bP_{t_0^*}) + \sqrt{2 K(\bP_{t_j^*}, \bP_{t_0^*})}\right]
 \end{align*}where $K(\bP_{t_j^*}, \bP_{t_0^*})$ denotes the Kullback-Leiber divergence between $\bP_{t_j^*}$ and  $\bP_{t_0^*}$. Since the noise is Gaussian and independent of $X$ in \eqref{eq:model} and $X$ is isotropic, we have $K(\bP_{t_j^*}, \bP_{t_0^*}) = N \norm{t_j^*-t_0^*}_2^2/(2\sigma^2) \leq N 9\rho^2/(8 \sigma^2)$. Hence, if $\rho\leq 2\sigma \sqrt{\log(1/\tau)/(96N)}$ and $\log(1/\tau)\geq 544/225$, one has 
 \begin{equation*}
 \bP_{t_j^*}\left[ \frac{d\bP_{t_0^*}}{d\bP_{t_j^*}}\geq \tau\right]\geq \frac{3}{4}.
 \end{equation*}

 Using the latter result together with \eqref{eq:main_adapt} in \eqref{eq:minimax_small_ball_main} for the values $\tau=1/(2d)$, we obtain that
 \begin{equation*}
 \bP_{t_0^*}\left[\norm{\hat t - t_0^*}_2\geq \rho/4\right]\geq\frac{1}{2}.
 \end{equation*}
 
 \endproof

\section{Minimax regularization function in the fixed design setup}
\label{fixed_design}
In this section, we consider the Gaussian regression model with fixed design. Our aim is to design a minimax regularization function in this setup. We are therefore given a deterministic $N\times d$ design matrix $\mathbb{X}$, whose $i$-th row vector is denoted by  $X_i\in\R^d$, and we observe $N$ noisy projections of a vector $t^*\in\R^d$,  $Y_i = \inr{X_i,t^*} + \xi_i, i=1, \ldots, N$, where the noises $\xi_1, \ldots, \xi_N$ are independent centered Gaussian variables with variances $\sigma^2$. The data $(Y_i, X_i)_{i=1}^N$ are used to construct estimators of $t^*$ and the only difference with the previous random design setup is that the $X_i$'s are deterministic vectors whereas so far they were random vectors. We will use most of the time the matrix form of the data: $Y=\bX t + \xi$ where $Y=(Y_i)_{i=1}^N$ and $\xi\sim\cN(0, \sigma^2 I_N)$ where $I_N$ is the $N\times N$ identity matrix. Note that the fixed design setup is usually considered in signal processing over a finite grid or in experiences where the statistician can chose in advance the design of an experience and then observed an output.

In order to design a minimax regularization function in this setup, we need to adapt the definitions introduced in Section~\ref{sec:introduction} to the fixed design case. We first start with the definition of a minimax rate of convergence over a subset of $\R^d$. We use the empirical (or normalized) Euclidean inner product and norm: for all $u, v\in\R^N$, 
\begin{equation*}
\inr{u, v}_{L_2^N} = \frac{1}{N}\sum_{i=1}^N u_i v_i \mbox{ and } \norm{u}_{L^2_N} = \sqrt{\frac{1}{N}\sum_{i=1}^N u_i^2}
\end{equation*}and the associated unit ball $B_{L^2_N} = \{u\in\R^N: \norm{u}_{L^2_N}\leq 1\}$.
 
\begin{Definition}\label{def:minimax fixed design}
Let $T\subset\R^d$,  $\mathbb{X}$ denote a (deterministic) $N\times d$ design matrix and $\xi$ be a centered Gaussian random vector in $\R^N$  with covariance matrix $\sigma^2 I_N$. For all $t^*\in T$, define the random vector $Y^{t^*}= \mathbb{X} t^* + \xi$ and denote by $\cY^T:=\{Y^{t^*}:t^*\in T\}$ the set of all such random vectors. 

Let $\hat t_N$ be a statistics from $\R^N$ to $\R^d$. Let $0<\delta_N<1$  and $\zeta_N>0$. We say that \textbf{$\hat t_N$ performs with accuracy $\zeta_N$ and confidence $1-\delta_N$ relative to the set of targets $\cY^T$}, if for all $t^* \in T$, with probability, w.r.t. to a vector  $Y$ distributed as $Y^{t^*}$, at least $1-\delta_N$,  $\norm{\bX(\hat t_N(Y)-t^*)}_{L^2_N}^2\leq \zeta_N$.

We say that \textbf{$\cR_N$ is a minimax rate of convergence over $T$ for the confidence $1-\delta_N$} if the two following hold:
\begin{enumerate}
	\item there exists a statistics $\hat t_N$ which performs with accuracy $\cR_N$ and confidence $1-\delta_N$ relative to the set of targets $\cY^T$
	\item there exists an absolute constant $g_0^\prime>0$ such that if $\tilde t_N$ is a statistics which performs with accuracy $\zeta_N$ and confidence $1-\delta_N$ relative to the set of targets $\cY^T$ then necessarily $\zeta_N \geq g_0^\prime \cR_N$.
\end{enumerate}
\end{Definition}
Note that we use the empirical $L_2^N$-metric $\norm{\bX \cdot}_{L^2_N}$ (to the square) with respect to the design $\bX$ as a measure of performances of estimators in Definition~\ref{def:minimax fixed design}. The reason we do so is that it is the natural counterpart to the random design case -- that is when $\bX$ is a standard Gaussian matrix then $\E\norm{\bX t}_{L^2_N}^2 = \norm{t}_2^2$ and the $\ell_2^d$-norm is the metric used to measure the performance of estimators in the random design setup -- and that it is the natural metric associated to the prediction of $Y$ problem given that if $R(t) = \E \norm{Y-\bX t}_{L^2_N}^2$ is the risk of $t$ for all $t\in\R^d$ then we have for any estimator $\hat t_N$, $R(\hat t_N)=R(t^*)+\norm{\bX(\hat t_N-t^*)}_{L^2_N}^2$. So that predicting $Y$ is the same problem as estimating $t^*$ with respect to the  empirical $\norm{\bX \cdot}_{L^2_N}$ metric.

Now, we adapt the definition of a minimax regularization function to the fixed design setup in the next definition. 

\begin{Definition}\label{def:optimality fixed design}
Let $\norm{\cdot}$ be a norm on $\R^d$, $T\subset \R^d$ and $0<\delta_N<1$.   Let us consider the following RERM for some given function $\Psi:\R_+\to \R$:
\begin{equation*}
\hat t \in\argmin_{t\in \R^d}\left(\norm{Y-\bX t}_{L^2_N}^2+\Psi(\norm{t})\right)
\end{equation*}constructed from a $N\times d$ deterministic matrix $\mathbb{X}$ and a random vector $Y=\mathbb{X}t^*+\xi$, with $\xi \sim \sigma \cN(0,I_N)$. We say that $\Psi$ is a \textbf{minimax regularization function for the norm $\norm{\cdot}$ over $T$ for the confidence $1-\delta_N$}, if there exists an absolute constant $g_1^\prime>0$ such that for all $t^*\in T$, the RERM $\hat t$ is such that with probability at least $1-\delta_N$, $\norm{\hat t -t^*}_{L^2_N}^2\leq g_1^\prime\cR_N$, where $\cR_N$ is the minimax rate of convergence over $\{t\in \R^d: \norm{t}\leq \norm{t^*}\}$.
\end{Definition}


The statistical bounds one can prove in the fixed design setup depend generally on the property of the design matrix $\mathbb{X}$. Many different assumptions have been introduced during the last two decades in high-dimensional statistics and we refer to \cite{van2009conditions} for some of them. In particular, norm preserving properties like the RIP or weaker assumption on the restricted eigenvalues like the REC or CC have played an important role in statistics (cf. \cite{MR2723472,deter_lasso,MR2533469, giraud2014introduction,MR2576316,MR2807761}). In this paper, we assume that $\bX$ satisfies the  ``Restricted Isometry Property''. It appears that this condition is  equivalent (up to constants) to the property satisfied by a standard Gaussian matrix as in Lemma~\ref{prop:event_Omega^*}.

\begin{Assumption}[RIP($s$)]\label{ass:rip}
If $N<d$ and $N/\log(ed/N)>1$ then we set $s=N/\log(ed/N)$. We assume that all $t$ in $\Sigma_{s}:= \{x \in \mathbb{R}^d: \norm{x}_0 \leq s \}$ is such that
\begin{equation}
\label{eq:isometry}
 \frac{1}{2}\norm{t}_2 \leq \norm{\mathbb{X}t}_{L^2_N} \leq \frac{3}{2}\norm{t}_2
\end{equation} where $\norm{x}_0$ is the size of the support of $x$.  If $N\geq d$, we assume that \eqref{eq:isometry} is satisfied for all $t\in\R^d$.
\end{Assumption}

Note that in the high-dimensional case, i.e. $d>N$, only the situation $N/\log(ed/N)>1$ is considered in Assumption~\ref{ass:rip} to avoid the ultra-high dimensional phenomena discovered in \cite{MR2879672}. RIP was introduced in~\cite{MR2243152} and it has been widely used and discussed (cf. for example \cite{davenport2010analysis}, \cite{baraniuk2008simple} or \cite{garg2009gradient}), in particular in the field of Compressed Sensing. From our perspective, we use this result for two reasons: 1) the minimax results over $\ell_1^d$ balls we need to develop for our proof of minimax regularization function has been obtained in the fixed design under this condition (or an equivalent one) in \cite{MR2816337};  2) the complexity parameter that we will be using in the fixed design setup have been computed under very general design matrix $\bX$ in \cite{pierre} but they were only proved to be optimal under the RIP assumption. Our main result covers more general design matrix than the one satisfying RIP but it turns out that those results do not allow to conclude on the minimax optimality of the associated regularization function beyond the RIP case; moreover and to our knowledge no sharp closed form are available for the computation of this regularization function beyond the RIP case.    

 
 \subsubsection{The multiplier process and its associated fixed point}
\label{sec:multi_proc_fixed}  
 Our analysis is based upon the study of the same regularized excess empirical risk quantity as in the random design section: for all $t\in\R^d$, 
 $$P_N\cL_t^\Psi:=\left(\norm{Y-\bX t}_{L^2_N}^2+ \Psi(t)\right)- \left(\norm{Y-\bX t^*}_{L^2_N}^2 + \Psi(t^*)\right).$$ 
We use the same quadratic / multiplier decomposition as in the random design case: for all $t\in\R^d$,  $P_N \cL_t^\Psi = P_N \cQ_{t-t^*} + P_N \cM_{t-t^*} + \cR_{t-t^*}$ where
 \begin{itemize}
\item  $P_N {\cal Q}_{t-t^*}:= \norm{\mathbb{X}(t-t^*)}_{L^2_N}^2$ is the quadratic part
\item $P_N {\cal M}_{t-t^*}:=2\inr{\xi, \bX(t^*-t)}$ is  the multiplier part
\item  ${\cal R}_{t,t^*}:=\Psi(t)-\Psi(t^*)$ is the regularization part.
\end{itemize}
 Contrary to the random design case, in the fixed design setup the only source of randomness is the Gaussian noise $\xi$, in particular the quadratic process is fully deterministic. Therefore, there is no need to define  a fixed point similar to $r_Q$ for a control on the quadratic process.  The only fixed point we introduce is a version of the previous multiplier fixed point $r_M(\cdot)$ adapted to the fixed design setup: for $\eta^\prime=  1/8$, let  
\begin{equation}\label{eq:fixed_point_multi_fixed}
r_\mathbb{X}(\rho):= \inf\left(r>0:\sigma \ell^*\left(\frac{\rho}{\sqrt{N}}\mathbb{X} B_1^d \cap r B_2^{N}\right)\leq \eta^\prime r^2 \sqrt{N} \right)
\end{equation} where $\mathbb{X}B_1^d:=\{\mathbb{X}t: t \in B_1^d\}$ and $B_2^N$ is the unit ball in $\ell_2^N$.

Define $r_0^\prime$ as the non-zero solution to the equation $\sigma \ell^*(rB_2^{N}\cap \im(\bX))=\eta^\prime r^2 \sqrt{N}$, where $\im(\bX)$ is the image  of $\bX$ in $\mathbb{R}^d$, i.e. 
\begin{equation*}
r_0^\prime = \sigma \ell^*(B_2^{N}\cap \im(\bX)) / (\eta^\prime\sqrt{N}) = \sigma \ell^*(B_2^{\rk(\bX)})/(\eta^\prime \sqrt{N}) = (\sigma/\eta^\prime) \sqrt{\rk (\bX)/N}
\end{equation*}
Let $\rho_0^\prime$ be the smallest $\rho$ such that $(\rho/\sqrt{N})\mathbb{X}  B_1^d$ contains $r_0^\prime B_2^N\cap \im(\bX)$. An argument  similar to the one used in the random design case shows that $\rho_0^\prime$ is the smallest radius such that for all $\rho\geq \rho_0^\prime$, $r_\bX(\rho) = r_\bX(\rho_0^\prime) = r_0^\prime$. Finally, we consider $K_0^\prime=\min \{k \in \mathbb{N}: 2^k \rho^* \geq 2\rho_0^\prime\}$.

The fixed point function $r_\bX(\cdot)$ depends on the Gaussian mean width of $\bX B_1^d$ intersected with $rB_2^N$  for various radii $r$. This quantity has been recently controlled in Proposition~2 in \cite{pierre}. 

\begin{Proposition}[Proposition~2 in \cite{pierre}]\label{prop:pierre} Let $\bX\in\R^{N \times d}$. Assume that the column vectors of $\bX$ are in $B_2^N$. Then, for all $r\geq0$,
\begin{equation*}
\ell^*(\bX B_1^d\cap r B_2^N)\leq \min\left(4 \sqrt{\log(8ed)}, 4 \sqrt{\log(8ed r^2)}, r\sqrt{\rk(\bX)}\right).
\end{equation*}
\end{Proposition}

It follows from some calculations (similar to the one used to obtain the closed form of $r_M(\cdot)$ in \eqref{eq:rM}) that there exists an absolute constant $C_\bX$ such that for all $\rho$, $r_\bX(\rho)\leq \overline{r_\bX}(\rho)$ with 
\begin{equation}\label{eq:rM_fixed}
 \overline{r_\bX}^2(\rho) = C \left\{
\begin{array}{cl}
\frac{\sigma^2 \rk(\bX)}{N} & \mbox{ if } \rho \geq \rho_0^\prime \\
\\
\rho \sigma\sqrt{\frac{1}{N}\log\Big(\frac{e\sigma d}{\rho\sqrt{N}}\Big)} & \mbox{ if } \sigma^2 \log d \leq \rho^2N\leq  \rho_0^{\prime2} N\\
\\
 \rho \sigma \sqrt{\frac{\log(ed)}{N}} & \mbox{ if } \rho^2N\leq \sigma^2 \log d.
\end{array}
\right.
\end{equation} We see that $\overline{r_\bX}(\rho)$ in \eqref{eq:rM_fixed} and $r_M(\rho)$ in \eqref{eq:rM} are very close. The only difference comes from the rank of $\bX$ and when $N\geq \zeta^\prime d$ and  $\rk(\bX)\sim d$, the two fixed points $r_M$ and $\overline{r_\bX}$ are equal up to absolute constants. Furthermore,  one can check that there are two absolute constants $0<C_1<C_2$ and $C_\bX=C_\bX(\rho, d, \sigma, N)$ such that $C_1\leq C_\bX\leq C_2$ and 
 \begin{equation}
 \label{eq: def rbarre}
  \overline{r_\bX}^2(\rho)= C_\bX \min \left(\frac{\sigma^2 \rk(\bX)}{N},\rho \sigma\sqrt{\frac{1}{N}\log\left(\frac{e\sigma d}{\rho\sqrt{N}}\right)}, \rho \sigma \sqrt{\frac{\log(ed)}{N}}\right).
 \end{equation} 

We will use $\rho\to\overline{r_\bX}^2(\rho)$ as a regularization function (up to an absolute constant). It will appear that  even though $\rho\to\overline{r_\bX}^2(\rho)$  is only an upper bound on  $\rho\to r_\bX(\rho)$, it will be an optimal minimax choice when the design matrix satisfies RIP.

 \subsubsection{Main result: a minimax regularization function}
 In this section, we prove that the regularization function $\Psi(\rho) = c_0^\prime \overline{r_\bX}^2(\rho)$ is a minimax regularization function for some well-chosen absolute constant $c_0^\prime$ when the design matrix $\bX$ satisfies RIP. To that end, we first need to know the minimax rate over $\ell_1^d$-balls in the fixed design setup. Such a result was obtained in \cite{MR2816337}. Let us now recall this result in our context (see (5.25) in Section~5.2.2 in \cite{MR2816337}).  

\begin{Proposition}[\cite{MR2816337}]\label{prop:sacha_philippe}
Let $\bX\in\R^{N \times d}$ be a matrix satisfying RIP($2s$). For all $\rho\geq0$, the minimax rate of convergence over $\rho B_1^d$ in the Gaussian linear model with fixed design $\bX$ in expectation is given by
\begin{equation}\label{eq:minimax_rate}
\min \left(\frac{\sigma^2 \rk(\bX)}{N},\rho \sigma\sqrt{\frac{1}{N}\log\left(\frac{e\sigma d}{\rho\sqrt{N}}\right)}, \rho^2\right).
\end{equation}
\end{Proposition}
The latter result holds in expectation whereas we are interested in deviation results. Even though the minimax rate of convergence in deviation over $\rho B_1^d$ in the Gaussian linear model under RIP has not been established, we believe that this rate of convergence in deviation is identical to the one given in Proposition~\ref{prop:sacha_philippe}. Note that a proof of this fact follows from Section~3 in \cite{LM13} for the minimax lower bound in deviation and from the quadratic / multiplier decomposition of the excess loss together with Lemma~\ref{bound on PN} below to show that the ERM over $\rho B_1^N$ achieves the minimax bound in deviation. We do not provide the proof here but we will use \eqref{eq:minimax_rate} as a benchmark result to construct a minimax regularization function.

\begin{Theorem}\label{theo:main fixed design}Let $\bX\in\R^{N\times d}$ be such that the column vectors of $\bX$ are in $B_2^N$. Consider the following regularization function:
 \begin{equation*}
\rho\geq 0\to \Psi(\rho) = c_0^\prime \overline{r_\bX}^2(\rho)
\end{equation*}
where $\overline{r_\bX}(\rho)$ is defined in~\eqref{eq: def rbarre}  and $c_0^\prime\geq 2$ is an absolute constant. Then there exist absolute constants $\kappa_1^\prime$, $\kappa_2^\prime$ and $\kappa_3^\prime$ such that for any $t^*\in\R^d$ the RERM $\hat t$ constructed from the data $Y=\bX t^*+\xi$:  
$$\hat t\in\argmin_{t\in \R^d}\left(\norm{Y-\bX t}_{L^2_N}^2+\Psi(\norm{t}_1)\right)$$  satisfies with probability greater than $1-\kappa_1^\prime\exp\left(-\kappa_2^\prime N \overline{r_\mathbb{X}}(\norm{t^*}_1)^2/\sigma^2\right)$,
 $$\norm{\bX(\hat t - t^*)}_{L_2^N}^2\leq  \kappa_3^\prime \min\left(\frac{\sigma^2 \rk(\mathbb{X})} {N}, \frac{\norm{t^*}_1 \sigma}{\sqrt{N}}\sqrt{\log\left(\frac{ed\sigma}{\norm{t^*}_1\sqrt{N}}\right)}, \norm{t^*}_1\sigma \sqrt{\frac{\log(ed)}{N}}, \right).$$ 
 \end{Theorem}
 
 Note that the probability estimate $1-\kappa_1^\prime\exp\left(-\kappa_2^\prime N \overline{r_\mathbb{X}}(\norm{t^*}_1)^2/\sigma^2\right)\geq 3/4$ only when $\norm{t^*}_1\geq  \Delta_0^\prime \sigma \sqrt{\log(ed)/N}$ for some absolute constant $\Delta_0^\prime$ large enough. Therefore, if \eqref{eq:minimax_rate} is indeed the minimax rate of convergence over $\rho B_1^d$ for the deviation $1-\delta_N=3/4$ under RIP then Theorem~\ref{theo:main fixed design} proves that $\rho\geq 0\to \Psi(\rho) = c_0^\prime \overline{r_\bX}^2(\rho)$ is a minimax regularization function for the $\ell_1^d$-norm over $\R^d\backslash \Delta_0^\prime \sigma \sqrt{\log(ed)/N} B_1^d$ for the constant confidence regime.

 \subsection{Proof of Theorem~\ref{theo:main fixed design}}
 The proof is split into a probabilist part used to identify a high probability event $\Omega^\prime_0$ on which the multiplier process is well controlled on the entire space $\R^d$ and a deterministic part where it is proved that, on the event $\Omega^\prime_0$, $P_N \cL_t^\Psi>0$ if $\norm{\bX(t-t^*)}_{L^2_N}^2\gtrsim \overline{r_\bX}(\norm{t^*}_1)^2$.

\subsubsection{Probabilistic control of the multiplier process}

The following lemma shows how the fixed point $r_\bX$ allows to control the multiplier process.

\begin{Lemma}
\label{bound on PN}
Let $\rho>0$ and take $\eta^\prime = 1/8$. Then, for all $r\geq r_\mathbb{X}(\rho)$, with probability greater than $1-\exp\left(-N r^2/(128\sigma^2)\right)$, for all $t \in t^*+ \rho B_1^d$,
\begin{equation*}
|P_N {\cal M}_{t-t^*}| \leq  \frac{1}{2}\max\left( r^2, \norm{\bX(t-t^*)}_{L^2_N}^2  \right).
\end{equation*}
\end{Lemma}

\beginproof 
Let $r\geq r_\mathbb{X}(\rho)$. We denote by $B_\bX^2$ the unit ball associated with the pseudo-metric $\norm{\bX\cdot}_{L^2_N}$ and $r B_\bX^2 = \left\{t \in \R^d: \norm{\bX t}_{L^2_N} \leq r\right\}$ its unit ball of radius $r$. We have
$$
\sup_{t\in t^*+\rho B_1^d\cap r B_{\bX}^2}|P_N {\cal M}_{t-t^*}|=
2\sup_{t\in t^*+\rho B_1^d\cap r B_{\bX}^2} \left|\inr{\mathbb{X}(t-t^*),\xi}_{L^2_N}\right|= \frac{2}{\sqrt{N}}\sup_{v \in V} \inr{v,\xi}$$ where $V= [(\rho/\sqrt{N}) \mathbb{X}B_1^d] \cap r B_2^{N}$. Then, it follows from Borell's concentration inequality (cf. \cite{MR1849347}) that for all $x>0$, with probability at least $1-\exp(-x^2/2)$, 
 $$ \sup_{v \in V} \inr{v, \xi} \leq  \E\sup_{v \in V} \inr{v, \xi} +  x \sigma(V)$$
where $\sigma(V):=\sup_{v \in V} \sqrt{\E \inr{v, \xi}^2} = \sup_{v\in V} \sigma \norm{v}_2 \leq \sigma r$. Moreover, given that $\xi\sim\cN(0, \sigma^2 I_N)$, we have 
\begin{equation*}
\E\sup_{v \in V} \inr{v, \xi} = \sigma \ell^*(V) = \sigma \ell^* \left(\frac{\rho}{\sqrt{N}}\mathbb{X} B_1^d \cap r B_2^{N}\right)\leq \eta^\prime r^2 \sqrt{N}
\end{equation*}where the last inequality follows from the definition of $r_\bX(\rho)$ and because $r\geq r_\bX(\rho)$. Gathering all the pieces together, it follows for $x= \eta^\prime r \sqrt{N} /\sigma$, that, with probability at least $1-\exp(-r^2 N/(128\sigma^2))$,
\begin{equation}\label{eq:event_multi_fixed}
\sup_{t\in t^*+\rho B_1^d\cap r B_{\bX}^2} \left| P_N \cM_{t-t^*}\right|\leq 4 \eta^\prime r^2 =  \frac{1}{2} r^2.
\end{equation}

Now, the proof follows from an homogeneity argument. Indeed, let us assume that \eqref{eq:event_multi_fixed} holds. Let $t \in t^*+\rho B_1^d$ be such that $\norm{\bX(t-t^*)}_{L^2_N}>r$ and define $t^\prime:=t^* + \alpha_t(t-t^*)$ where $\alpha_t = r/\norm{\bX(t-t^*)}_{L^2_N}$. Note that $\alpha_t<1$ and  $t^\prime \in t^*+\rho B_1^d\cap r B_{\bX}^2$. Hence, it follows from \eqref{eq:event_multi_fixed} that $|P_N {\cal M}_{t^\prime-t^*}| \leq  r^2/2$ and so $|P_N {\cal M}_{t-t^*}| = |P_N {\cal M}_{t^\prime-t^*}| / \alpha_t \leq  r^2/(2\alpha_t)\leq \norm{\bX(\hat t-t^*)}_{L^2_N}^2/2$.
\endproof



\subsubsection{Deterministic part of the proof}
We start with two lemmas on the growth behavior of $\overline{r_\bX}(\cdot)$. Their proofs are almost identical to the one of Lemmas~\ref{rate of growth} and~\ref{reverse Lemma} and are therefore omitted. 

\begin{Lemma}
\label{rate of growth, fixed design} Let $\phi^\prime= 4$. If $\phi^\prime\rho \leq \rho_0^\prime\,  \min(1,\eta^\prime)$, then for any $\rho^\prime\geq \phi^\prime  \rho$,
$\overline{r_\bX}^2\left(\rho^\prime\right)>2\overline{r_\bX}^2(\rho).$
\end{Lemma}

\begin{Lemma}
\label{reverse Lemma, fixed design} Let $\nu>0$. 
If $\nu \geq 1$ then $\overline{r_\bX}(\nu \rho)\leq \sqrt{\nu}\overline{r_\bX}(\rho)$. If $\nu \leq 1$ then $\overline{r_\bX}(\nu \rho)\geq \sqrt{\nu}\overline{r_\bX}(\rho)$.
\end{Lemma}

To prove Theorem~\ref{theo:main fixed design}, we use the same argument as in the proof of Theorem~\ref{theo:main} for the random design. Let $\rho^*=10 \norm{t^*}_1/\eta^\prime$ 
and split $\R^d$ into three zones:
\begin{itemize}
\item the ``central zone" $t^*+\rho^*B_1^d$,  
\item the intermediate ``peeling zone": $\{t\in\R^d: \rho^*<\norm{t-t^*}_1\leq 2^{K_0^\prime}\rho^*\}$  -- to be considered only when $K_0^\prime \geq 1$. This part of $\R^d$ is itself partitioned into $K_0^\prime$ shelves: for $k=1, \ldots, K_0^\prime$, $\{t\in\R^d: 2^{k-1} \rho^*<\norm{t-t^*}_1\leq 2^{k} \rho^*\}$, 
\item the ``exterior zone": $\{t\in\R^d: \norm{t-t^*}_1> 2^{K_0^\prime}\rho^*\}$ on which $r_\bX$ is constant equal to $r_0^\prime$.
\end{itemize}
For all $ k = 0,\ldots, K_0^\prime$, we denote by $A^\prime_k$ the event on which for all $t \in t^*+ 2^k\rho^* B_1^d$,
\begin{equation*}
|P_N {\cal M}_{t-t^*}| \leq  \frac{1}{2}\max\left( \overline{r_\mathbb{X}}(2^k \rho^*)^2, \norm{\bX(t-t^*)}_{L^2_N}^2  \right).
\end{equation*} We consider the event $\Omega_0^\prime =  A^\prime_0\cap\cdots\cap A^\prime_{K_0^\prime}$. It follows from Lemma~\ref{bound on PN} and an argument similar to the one in Lemma~\ref{pseudo-geometric} that for some absolute constants $\kappa_1^\prime$, $\kappa_2^\prime$ and $\kappa_4^\prime$, 
$$\bP\left[\Omega_0^\prime\right] \geq 1- \kappa_1^\prime \exp\big(-\kappa_2^\prime N \overline{r_\bX}(\norm{t^*}_1)^2/\sigma^2\big)$$ as long as $\norm{t^*}_1\geq \kappa_4^\prime \sigma/\sqrt{N}$ (which is the case when $\norm{t^*}_1\geq \Delta_0^\prime \sigma \sqrt{\log(ed)/N}$ for $\Delta_0^\prime\geq \kappa_4^\prime$). 

Let us now assume for the remaining of the proof that $\Omega_0^\prime$ holds. Note that unlike in the random design case, there is no event such as $\Omega^*$ in $\Omega_0^\prime$ on which the quadratic process is controlled, because, in the deterministic design case this process is deterministic.

Our strategy is to show that $\hat t$ belongs to the ``central zone''. To that end it is enough to prove that $P_N \cL_t^\Psi>0$ for every $t\in\R^d$ such that $\norm{t-t^*}_1> \rho^*$ because by definition $P_N \cL_{\hat t}^\Psi\leq0$. 

Let $t$ be in the intermediate peeling zone (which can happen only if $\rho^*\leq 2 \rho_0^\prime$), say in the $k$-th shell for some $k\in\{0, \ldots, K_0^\prime\}$: $2^{k-1} \rho^*< \norm{t-t^*}_1 \leq 2^{k} \rho^*$.  In particular $\norm{t}_1 > \norm{t^*}_1$ and ${\cal R}_{t,t^*}\geq 0$. Therefore, if $\norm{\bX(t-t^*)}_{L^2_N} \geq \overline{r_\bX}(2^{k} \rho^*)$ then by Lemma~\ref{bound on PN}, $ |P_N {\cal M}_{t-t^*}|\leq \norm{\bX(t-t^*)}_{L^2_N}^2=P_N {\cal Q}_{t-t^*}$ and so $P_N\cL_t^\Psi>0$. Now, if $\norm{\bX(t-t^*)}_{L^2_N} \leq \overline{r_\bX}(2^{k} \rho^*)$ then by Lemmas~\ref{bound on PN} and~\ref{reverse Lemma, fixed design}, for $c_0^\prime \geq 2$,
$$ |P_N {\cal M}_{t-t^*}| \leq \frac{1}{2} \overline{r_\bX}(2^{k}\rho^*)^2  < \frac{c_0^\prime}{2} \overline{r_\bX}^2(\norm{t-t^*}_1)$$
and since $\norm{t}_1 \geq \norm{t-t^*}_1-\norm{t^*}_1\geq 4\norm{t^*}_1$, and $4\norm{t^*}_1 \leq \eta^\prime\rho^*/2 \leq \eta^\prime\rho_0^\prime$, by Lemma~\ref{rate of growth, fixed design} , one has $c_0^\prime \overline{r_\bX}^2(\norm{t}_1)\geq 2c_0^\prime \overline{r_\bX}^2(\norm{t^*}_1)$. As a consequence, ${\cal R}_{t,t^*} \geq c_0^\prime \overline{r_\bX}^2(\norm{t}_1)/2> |P_N {\cal M}_{t-t^*}|$ and so $P_N\cL_t^\Psi>0$.

Let us now tackle the exterior zone in both cases $\rho^* \leq 2 \rho_0^\prime$ and $\rho^* > 2 \rho_0^\prime$ .  Let $t\in\R^d$ be such that  $\norm{t-t^*}_1 >  2^{K_0^\prime} \rho^*$. We have ${\cal R}_{ t,t^*} \geq 0$ because $\Psi(\norm{t}_1) = c_0^\prime \overline{r_\bX}^2(\norm{t}_1)=c_0^\prime  r_0^{\prime 2} \geq c_0^\prime \overline{r_\bX}^2(\norm{t^*}_1) = \Psi(\norm{t^*}_1)$. Let $t^\prime = t^*+\alpha_t(t-t^*)$ for some $0<\alpha_t<1$ be such that $\norm{t^\prime - t^*}_1 = 2^{K_0^\prime} \rho^*$. By definition of $K_0^\prime$ and $\rho_0^\prime$, we have $\norm{\bX(t^\prime - t^*)}_{L^2_N} \geq  r_0^\prime$. Therefore,  since $A_{K_0^\prime}\subset \Omega_0^\prime$, we have $|P_N {\cal M}_{t^\prime-t^*}|\leq (1/2)\norm{\bX(t^\prime-t^*)}_{L^2_N}^2$ which implies that $P_N \cQ_{t^\prime - t^*} + P_N\cM_{t^\prime - t^*}>0$ and therefore by an homogeneity argument that $P_N \cQ_{t - t^*} + P_N\cM_{t - t^*}>0$. Finally, given that ${\cal R}_{\hat t,t^*} \geq 0$ we conclude that $P_N \cL_{t}^\Psi>0$. 

This proves that $\hat t$ lies in the central zone in both cases $\rho^*\leq 2\rho_0^\prime$ and $\rho^*> 2\rho_0^\prime$.  But, now given that $A_0^\prime\subset \Omega_0^\prime$, we have
 $$|P_N {\cal M}_{\hat t-t^*}|\leq \frac{1}{2}\max\left( \overline{r_\bX}(\rho^*)^2, \norm{\bX(\hat t-t^*)}_{L^2_N}^2  \right).$$
 If $\norm{\bX(\hat t-t^*)}_{L^2_N}^2\leq \overline{r_\bX}(\rho^*)^2$ the proof is over and otherwise $|P_N {\cal M}_{\hat t-t^*}|\leq (1/2) \norm{\bX(\hat t-t^*)}_{L^2_N}^2$ which implies that $\norm{\hat t - t^*}_{L^2_N}^2\leq 2 \Psi(\norm{t^*}_1)\leq 2 c_0^\prime \overline{r_\bX}(\rho^*)^2$ because
\begin{equation*}
0\geq P_N\cL_{\hat t}^\Psi\geq \frac{1}{2}\norm{\hat t - t^*}_{L^2_N}^2 + \Psi(\norm{\hat t}_1) -  \Psi(\norm{t^*}_1). 
\end{equation*}This proves, on the event $\Omega_0^\prime$, that $\norm{\hat t - t^*}_{L^2_N}^2\leq 2 c_0^\prime \overline{r_\bX}(\rho^*)^2$.
\endproof

\begin{footnotesize}
\bibliographystyle{plain}
\bibliography{biblio}
\end{footnotesize}

\end{document}